\documentclass[1p,sort&compress]{elsarticle}

\usepackage{amssymb}
\usepackage{amsmath, amsthm}
\usepackage{epsfig}
\usepackage{graphicx}

\makeatletter
\@addtoreset{equation}{section}
\makeatother

\newtheorem{thm}{Theorem}[section]
\newtheorem{lem}[thm]{Lemma}
\newtheorem{corollary}[thm]{Corollary}
\newtheorem{conjecture}[thm]{Conjecture}
\theoremstyle{definition}
\newtheorem{definition}[thm]{Definition}

\newproof{pf}{Proof}
\newdefinition{rmk}[thm]{Remark}

\begin{document}
\begin{frontmatter}
\title{On sizes of complete arcs in $PG(2,q)$}

\author[per]{Daniele Bartoli}
\ead{daniele.bartoli@dmi.unipg.it}
\author[mosc]{Alexander A. Davydov}
\ead{adav@iitp.ru }
\author[per]{Giorgio Faina}
\ead{faina@dmi.unipg.it}
\author[per]{Stefano Marcugini}
\ead{gino@dmi.unipg.it}
\author[per]{Fernanda Pambianco\corref{cor1}}
\ead{fernanda@dmi.unipg.it}
\address[per]{Dipartimento di Matematica e Informatica, Universit\`{a}
degli Studi di Perugia, Via Vanvitelli~1,\\ Perugia, 06123,
Italy}
\address[mosc]{Institute for Information Transmission Problems,
Russian Academy of Sciences,\\ Bol'shoi Karetnyi per. 19,
GSP-4, Moscow, 127994, Russia}
\cortext[cor1]{Corresponding author~phone +39(075)5855006 ~fax
+39(075)5855024
}
\begin{abstract}
New upper bounds on the smallest size $t_{2}(2,q)$ of a
complete arc in the projective plane $PG(2,q)$ are obtained for
$853\leq q\leq 4561$ and $q\in T_{1}\cup T_{2}$ where $T_{1}=
\{173,181,193,229,\linebreak
243,257,271,277,293,343,373,409,443,449,457,461,463,467,479,
487,491,499,529,\linebreak 563,
569,571,577,587,593,599,601,607,613,617,619,631,641,661,
673,677, 683,691,\linebreak 709\}$,
 $ T_{2}=
 \{4597,4703,4723,4733,4789,4799,4813,4831,5003,5347,5641,5843,6011,\\
 8192\}$.
 From these new bounds it follows that for
$q\leq 2593$ and $q=2693,2753$, the relation
$t_{2}(2,q)<4.5\sqrt{q}\smallskip $ holds. Also, for $q\leq
4561$ we have $t_{2}(2,q)<4.75\sqrt{q}.$ It is showed that for
$23\leq q\leq 4561$ and $q\in T_{2}\cup
\{2^{14},2^{15},2^{18}\}$, the inequality
$t_{2}(2,q)<\sqrt{q}\ln ^{0.75}q$ is true. Moreover, the
results obtained allow us to conjecture that this estimate
holds for all $q\geq 23.$ The new upper bounds are obtained by
finding new small complete arcs with the help of a computer
search using randomized greedy algorithms. Also new
constructions of complete arcs are proposed. These
constructions form families of $k$-arcs in $PG(2,q)$ containing
arcs of all sizes $k$ in a region $k_{min}\le k \le k_{max}$
where $k_{min}$ is of order $\frac{1}{3}q$ or $\frac{1}{4}q$
while $k_{max}$ has order $\frac{1}{2}q$. The completeness of
the arcs obtained by the new constructions is proved for $q\le
1367$ and $2003\le q\le2063$. There is reason to suppose that
the arcs are complete for all $q>1367$. New sizes of complete
arcs in $PG(2,q)$ are presented for $169\leq q\leq 349$ and
$q=1013,2003.$
\end{abstract}


\end{frontmatter}

\section{Introduction}

Let $PG(2,q)$ be the projective plane over the Galois field
$F_{q}$. An $k$-arc is a set of $k$ points no three of which
are collinear. An $k$-arc is called complete if it is not
contained in an $(k+1)$-arc of $PG(2,q)$. For an introduction
in projective geometries over finite fields, see
\cite{HirsBook,SegreLeGeom,SegreIntrodGalGeom}.

In \cite{HirsSt-old,HirsStor-2001} the close relationship
between the theory of $k$-arcs, coding theory and mathematical
statistics is presented. In particular, a complete arc in a
plane $PG(2,q),$ points of which are treated as 3-dimensional
$q$-ary columns, defines a parity check matrix of a $q$-ary
linear code with codimension 3, Hamming distance 4, and
covering radius 2. Arcs can be interpreted as linear maximum
distance separable (MDS) codes \cite{szoT93,thaJ92d} and they
are related to optimal coverings arrays \cite{Hartman-Haskin}
and to superregular matrices \cite{Keri}.

One of the main problems in the study of projective planes, which is also of
interest in coding theory, is finding of the spectrum of possible sizes of
complete arcs.

A great part of this work is devoted to upper bounds on
$t_{2}(2,q)$, the smallest size of a complete arc in $PG(2,q)$.
Also we propose new constructions of complete arcs in $PG(2,q)$
and consider the spectrum of their possible sizes.

Surveys of results on the sizes of plane complete arcs, methods
of their construction and comprehension of the relating
properties can be found in
\cite{DFMP-JG2005,DFMP-Pampor,DFMP-JG2009,FP,HirsSurvey83,HirsBook,HirsSt-old,HirsStor-2001,
LombRad,pelG77,pel83,pelG92,pel93,SegreLeGeom,SegreIntrodGalGeom,SZ,segB62b,szoT87a,szoT89survey,szoT93}.
In particular, as it is noted in
\cite{HirsSt-old,HirsStor-2001}, the following idea of Segre
\cite{segB62b} and Lombardo-Radice \cite{LombRad} is fruitful:
the points of the arc are chosen, with some exceptions, among
the points of a conic or a cubic curve. We use this idea for
constructions of complete arcs and for finding the spectrum of
arc sizes, see Sections~\ref {sec2_computSearch},
\ref{sec5_constr}, and \ref{sec6_spectrum}.

The maximum size $m_{2}(2,q)$ of a complete arc in $PG(2,q)$ is well known.
It holds that
\begin{equation*}
m_{2}(2,q)=\left\{
\begin{array}{ll}
q+1 & \text{if }q\text{ odd} \\
q+2 & \text{if }q\text{ even}
\end{array}
\right. .
\end{equation*}

On the other hand, finding an estimation of the minimum size $t_{2}(2,q)$ is
a hard open problem.

Problems connected with small complete plane arcs are
considered in \cite {DFMP-JG2005,DFMP-Pampor,DFMP-JG2009,
DGMP-Innov,DGMP-JCD,FainaGiul,FPDM,
Giul2000,Giul2007affin,MaPa-private,Giul2007even,GiulUghi,
HirsBook,KV,LisMarcPamb2008, MMP-q29,MaMiP32,Ost,
Polv,SegreLeGeom,SZ,segB62b,szoT87a,szoT89survey,szoT93}, see also the
references therein.

We denote the aggregates of $q$ values:
\begin{eqnarray*}
T_{1}
&=&\{173,181,193,229,243,257,271,277,293,343,373,409,443,449,457,461,463, \\
&&467,479,487,491,499,529,563,569,571,577,587,593,599,601,607,613,617,\\
&&619,631,641,661,673,677,683,691,709\};  \\
T_{2} &=&\{4597,4703,4723,4733,4789,4799,4813,4831,5003,5347,5641,5843,6011,\\
&&8192\};   \\
T_{3} &=&\{2^{14},2^{15},2^{18}\};\\
Q&=&\{961,1024,1369,1681,2401\}=\{31^{2},2^{10},37^{2},41^{2},7^{4}\};\\
N&=&\{601,5^{4},661,3^{6},29^{2},31^{2},2^{10},37^{2},41^{2},7^{4},2^{18}\}.
\end{eqnarray*}

The exact values of $t_{2}(2,q)$ are known only for $q\leq 32,$
see \cite{MMP-q29} and recent work \cite{MaMiP32} where the
equalities $ t_{2}(2,31)=t_{2}(2,32)=14$ are proven. Also,
there are the following lower bounds (see
\cite{Polv,SegreLeGeom}):
\begin{equation*}
t_{2}(2,q)>\left\{
\begin{array}{ll}
\sqrt{2q}+1 & \text{for any }q \\
\sqrt{3q}+\frac{1}{2} & \text{for }q=p^{h},\text{ }p\text{ prime, }h=1,2,3
\end{array}
\right. .
\end{equation*}

Let $\overline{t}_{2}(2,q)$ be the smallest \emph{known }size
of a complete arc in $PG(2,q)$. For $q\leq 841$ the values of
$\overline{t}_{2}(2,q)$ (up to June 2009) are collected in
\cite[ Tab.\thinspace 1]{DFMP-JG2009} whence it follows that
$\overline{t} _{2}(2,q)<4\sqrt{q}$ $\mbox{ for }$ $q\leq 841$ .
In \cite{Giul2000}, see also \cite{GiulUghi}, complete
$(4\sqrt{q}-4)$ -arcs are obtained for $q=p^{2} $ odd, $q\leq
1681$ or $q=2401$. In \cite {DGMP-JCD,FainaGiul} complete
$(4\sqrt{q}-4)$-arcs are obtained for even $q=64,256,1024.$ By
the results above-cited it holds that
\begin{equation}
t_{2}(2,q)<4\sqrt{q}\quad \mbox{ for }\text{ }2\leq q\leq 841,\text{ }
q\in Q.  \label{eq1_<4sqroot(q)}
\end{equation}

For even $q=2^{h}$, $11\leq h\leq 15$, the smallest known sizes
of complete $ n$ -arcs in $PG(2,q)$ are obtained in
\cite{DGMP-JCD}, see also \cite[ p.\thinspace 35]{DFMP-JG2009}.
They are as follows: $ \overline{t}_{2}(2,2^{11})=201,\text{
}\overline{t}_{2}(2,2^{12})=307,\text{
}\overline{t}_{2}(2,2^{13})=461,\text{
}\overline{t}_{2}(2,2^{14})=665,\text{
}\overline{t}_{2}(2,2^{15})=993.$ Also,
${(6\sqrt{q}-6)}${-arcs} in $PG(2,q),$ $q=4^{2h+1},$ are
constructed in \cite{DGMP-Innov}; for $h\leq 4$ it is proved
that they are complete. It gives a complete 3066-arc in
$PG(2,2^{18}).$

Let $t(\mathcal{P}_{q})$ be the size of the smallest complete arc in any
(not necessarily Galois) projective plane $\mathcal{P}_{q}$ of order $q$. In
\cite{KV}, for \emph{sufficiently large} $q$, the following result is proved
(we give it in the form of \cite[Tab.\thinspace 2.6]{HirsStor-2001}):
\begin{equation}
t(\mathcal{P}_{q})\leq d\sqrt{q}\log ^{c}q,\text{ }c\leq 300,
\label{eq1_KimVu_c=300}
\end{equation}
where $c$ and $d$ are constants independent of $q$ (i.e. universal
constants). The logarithm basis is not noted as the estimate is asymptotic.

In this work, by computer search using randomized greedy
algorithms (see Section \ref{sec2_computSearch}), new small
complete arcs in $PG(2,q)$ are obtained for all $q\in T_{1}\cup
T_{2}$ and for $853\leq q\leq 4561.$  For $ q=601,661,$ the new
complete arcs arose from a theoretical study on orbits of
subgroups, helped by computer \cite{MaPa-private}. From the
sizes of the new arcs, with the use of (\ref {eq1_<4sqroot(q)})
and \cite[Tab.\thinspace 1] {DFMP-JG2009}, \cite{DGMP-Innov},
the following theorems result (see also Theorems
\ref{Th3_4.5sqroot(q)}, \ref{Th3_5sqroot(q)}, and
\ref{th4_ln0.75} for more details).

\begin{thm}
\label{th1_<4.5_4.62} In $PG(2,q),$ the following holds.
\begin{eqnarray*}
t_{2}(2,q) &<&4.5\sqrt{q}~~\, \mbox{ for }q\leq
2593,\,\,q=2693,2753.\\
t_{2}(2,q) &<&4.75\sqrt{q}\,\,\mbox{ for }q\leq 4561.  \\
t_{2}(2,q) &<&4.98\sqrt{q}\,\,\mbox{ for }q\in T_{2}.
\end{eqnarray*}
\end{thm}

\begin{thm}
\label{th1_ln0.75} In $PG(2,q),$
\begin{equation}
t_{2}(2,q)<\sqrt{q}\ln ^{0.75}q\quad \mbox{ for }23\leq q\leq 4561,\text{ }
q\in T_{2}\cup T_{3}.  \label{eq1_ln0.75}
\end{equation}
\end{thm}

\noindent Moreover, the study of the values of
$\overline{t}_{2}(2,q)$ allows us to conjecture that the
estimate (\ref{eq1_ln0.75}) holds for all $ q\geq 23$ and the
last estimate of Theorem \ref{th1_<4.5_4.62} is right for all $
q\leq 8192$.

\begin{conjecture}
\label{conj1_ln0.75} In $PG(2,q),$
\begin{eqnarray}
t_{2}(2,q)&<&\sqrt{q}\ln ^{0.75}q\quad \mbox{ for }q\geq 23.
\label{eq1_conjecture_ln^0.75}\\
t_{2}(2,q) &<&5\sqrt{q}\hspace{1cm}\quad\mbox{ for }q\le8192.\notag
\end{eqnarray}
\end{conjecture}

Regarding the spectrum of complete arc sizes, we note (going
after \cite[ p.\thinspace 209]{HirsStor-2001}) that in
literature, complete arcs have constructed with sizes
approximately $\frac{1}{2}q$ (see
\cite{Belgium,DFMP-JG2009,DMP-ACCT2010,Giordano,HirsBook,KorchSon,
KorchSon2010,LombRad,pelG77,pelG92,pel93,SegreIntrodGalGeom}),
$\frac{1}{3}q$ (see \cite{abaV83,korG83a,SZ,szoT87a,Voloch90}),
$\frac{1}{4}q$ (see \cite{korG83a,szoT89survey}), $2q^{0.9}$
(see \cite{SZ} where such the arcs are constructed for
$q>7^{10}$). In particular, for even $ q\geq 8$ there is a
complete $\frac{1}{2}(q+4)$-arc \cite{HirsBook}. Important
results on the spectrum of complete arc sizes are collected in
\cite[Th.\thinspace 2.6]{HirsStor-2001} where it is noted, for
example, that in $PG(2,q)$ with $q$ square there exists a
complete $(q-\sqrt{q}+1)$-arc. In \cite{KorchPace}, \emph{large
complete arcs} in $PG(2,q^{n})$ are defined and new infinite families of the such arcs are
constructed.

Much attention is given to $\frac{1}{2}(q+5)$-arcs and
$\frac{1}{2}(q+7)$-arcs sharing $\frac{1}{2}(q+3)$ points with
a conic for $q$ odd \cite {Belgium,DFMP-Pampor,DFMP-JG2009},
\cite[Rem.\thinspace 3] {DMP-DESI},
\cite{DMP-ACCT2010,Giordano,KorchSon,KorchSon2010,pelG92}. It
is proved that for all odd $q$ there is a complete
$\frac{1}{2}(q+5)$-arc \cite{KorchSon2010}, see also
\cite{Belgium}. Also, a complete $\frac{1}{2}(q+7)$-arc exists
at least for the following \emph{odd }$q:$
\begin{eqnarray}
25 &\leq &q\leq 167\text{ \cite[Tab.\thinspace 2]{DFMP-JG2005},\cite[
Sec.\thinspace 2, Tab.\thinspace 2]{DFMP-JG2009},\cite[Tab.\thinspace 2.4]
{FP},\cite[Introduction]{KorchSon2010};}  \notag \\
q &\equiv &2\text{ }(\bmod~3),\text{ }q\leq 4523\text{ \cite{DFMP-JG2009,DMP-ACCT2010};}  \notag \\
q &\equiv &1\text{ }(\bmod~4),\text{ }q\leq 337\text{ \cite{Giordano};}
\label{eq1+(q+7)/2} \\
q &=&2bt-1,\text{ }t\text{ odd prime, }b=1,2\text{ \cite[Introduction]
{KorchSon2010};}  \notag \\
q &\equiv &3\text{ }(\bmod~4)\text{ and condition of \cite[Cor.\thinspace
4.16(i)]{KorchSon2010} holds, see \cite[Introduction]{KorchSon2010};}  \notag
\\
q^{2} &\equiv &1\text{ }(\bmod~16)\text{ and condition of \cite[
Cor.\thinspace 4.17]{KorchSon2010} holds, see \cite[Introduction]
{KorchSon2010}.}  \notag
\end{eqnarray}

For $q\leq 167$ the known sizes of complete arcs in $PG(2,q)$
are collected in \cite[Tab.\thinspace 2]{DFMP-JG2005},
\cite[Tab.\thinspace 2]{DFMP-JG2009}, \cite[Tab.\thinspace
2.4]{FP}.

In this work new Constructions A, B, C of complete arcs are
proposed, see Section~\ref{sec5_constr}. These constructions
form families of complete $k$-arcs in $PG(2,q)$ containing arcs
of all sizes $k$ in a region $k_{min}\le k \le k_{max}$ where
$k_{min}$ is of order $\frac{1}{3}q$ or $\frac{1}{4}q$ while
$k_{max}$ has order $\frac{1}{2}q$. The
completeness of the arcs obtained by the new constructions
is proved for $q\le 1367$ and $2003\le q\le2063$. Moreover, there is reason to suppose
that the arcs are complete for all $q>1367$. From Theorems
\ref{th5.1_result},
 \ref{th5.2_result}, \ref{th5.3_result} the following
theorem results:
\begin{thm}\label{th1_construc}
Constructions A,B, and C of Section \emph{\ref{sec5_constr}}
form families of complete $k$-arcs in $PG(2,q)$ containing arcs
of all sizes $k$ in the following regions:
\begin{eqnarray*}
\mbox{\emph{(i)}}~\mbox{Construction A}&:&\left\lfloor \frac{q+8}{3}\right\rfloor \leq k\leq \frac{q+5}{2},\,
q \mbox{ \emph{prime}},\\
&&109\leq q\leq 1367,\,2003\le q\le2063,\,q=73,97,101,103.\\
\mbox{\emph{(ii)}}~\mbox{Construction B}&:&
\left\lfloor \frac{q+8}{3}\right\rfloor \leq k\leq
\left\lfloor \frac{q+4}{2}\right\rfloor,q\not\equiv 3~(\bmod~4)
\mbox{ is  a \emph{prime power,}}\\
&&128\leq q\leq 1367,\,\,2003\le q\le2063,q=89,109,113,121.\\
\mbox{\emph{(iii)}}~\mbox{Construction C}&:&~
\frac{q+13}{4} \leq k\leq
 \frac{q+5}{2},\,~~~q\equiv 3~(\bmod~4)
\mbox{ is a \emph{prime power,}}\\
&&347\leq q\leq 1367,\,2003\le q\le2063,\\
&&q=199,227,239,243,251,263,271,283,307,311,331.
\end{eqnarray*}
\end{thm}
For the given $q$, in order to show that arcs obtained by a
construction are complete we should calculate by computer some
special value, say $\overline{L}_{q}$ (see Definitions
\ref{def5.1_overlineH}, \ref{def5.2_Gq}, \ref{def5.3_Jq}) and
check if $\overline{L}_{q}\le R_{q}$ where
$R_{q}=\lfloor\frac{1}{3}(q-1)\rfloor$ for Constructions A, B
and $R_{q}=\frac{1}{4}(q-3)$  for Construction C. The
calculations are relatively simple. Moreover, for $q\le1367$
and $2003\le q\le2063$ it holds that
$\overline{L}_{q}<\sqrt{q}\ln q$ and the difference
$R_{q}-\overline{L}_{q}$ has a tendency to increasing when $q$
grows, see Theorems \ref{th5.1_byTable1}(iii),(iv),
\ref{th5.2_byTable6}(iii),(iv), \ref{th5.3_byComput}(iii)(iv).
It allows us to conjecture the following, cf. Conjecture
\ref{conj5_constr}.
\begin{conjecture}
\label{conj1_construc} The assertions of Theorem
\emph{\ref{th1_construc}} hold also for all $q>1367$.
\end{conjecture}

A $k$-arc of Constructions A and B contains $k-2$ points in
common with a conic and two points lying on a tangent to the
conic (Construction A) or on a bisecant of the conic
(Construction B). A $k$-arc of Construction C contains $k-3$
points in common with a conic, two points lying on a bisecant
and one point on a tangent. Note that in \cite{abaV83},
$k$-arcs containing $k-2$ points of a hyperoval (the nucleus
among them) and two points on its bisecant are constructed.
Also, in the space $PG(3,q)$, $k$-caps with $k-2$ points in
common with a quadric are considered in
\cite{FP-JG1996,Pam-1999,PamUg-2000}. In
\cite{FP-JG1996,Pam-1999}, the cap contains two points on a
tangent to the quadric while in~\cite{PamUg-2000} two points
lie on an external line.

The complete arcs of Constructions A, B, C can be used as
starting objects for inductive constructions of \cite{DGMP-JCD} and
\cite{KorchPace}, see Remark \ref{rem5_starting-objects}. In that
way, using results of \cite{DGMP-JCD} together with Constructions A, B, C,
one can generate infinite sets of families of complete caps in
projective spaces $PG(v,2^{n})$ of growing dimensions $v$. Also, infinite families of large
complete arcs in $PG(2,q^{n})$ with growing $n$ can be obtained by constructions of
\cite{KorchPace} using arcs of Constructions A, B, C.

In this work, using Constructions A, B, C and randomized
greedy algorithms, new
complete arcs in $PG(2,q)$ are obtained for $ 169\leq q\leq
349$ and $q=1013,2003.$

In Section \ref{sec2_computSearch} we describe the greedy
algorithms used for obtaining new arcs. In Section~\ref{sec3_tables}
 we collect the known and new upper bounds on
$t_{2}(2,q)$ for $ q\leq 4561$ and $q\in T_{2}\cup T_{3}$. The
bounds are represented by tables, where values of
$\overline{t}_{2}(2,q)$ are written, and by the corresponding
relations. In Section \ref{sec4_observations} we give the upper
bounds on $ t_{2}(2,q) $ in the form of (\ref{eq1_ln0.75}) and
substantiate Conjecture \ref {conj1_ln0.75}. In Section
\ref{sec5_constr} new Constructions A, B, C of complete arcs
are described. Finally, in Section \ref{sec6_spectrum} we
present new sizes of complete arcs in $PG(2,q)$ with $169\leq
q\leq 349$ and $q=1013,2003$.

Some of the results of this work were briefly presented without
proofs in \cite{DFMP-ACCT2010}, see also~\cite{DFMP-ArXiv}.

\section{An approach to computer search\label{sec2_computSearch}}

In this paper for computer search we use the randomized greedy algorithms
\cite[Sec.\thinspace 2]{DFMP-JG2005},\cite[Sec.\thinspace 2]{DMP-JG2004}
that are convenient for relatively large $q$ and for obtaining examples of
different sizes of complete arcs. At every step an algorithm minimizes or
maximizes an objective function $f$ but some steps are executed in a random
manner. The number of these steps and their ordinal numbers have been taken
intuitively. Also, if the same extremum of $f$ can be obtained in distinct
ways, one way is chosen randomly.

We begin to construct a complete arc by using a starting set of
points $ S_{0} $. At the $i$-th step one point is added to the
set and we obtain a point set $S_{i}$. As the value of the
objective function $f$ we consider the number of points in
$PG(2,q)$ that lie on bisecants of the set obtained. For small
arcs we look for the maximum of the objective function $f.$ For
the spectrum of arc sizes we use both the maximum and the
minimum of $f$.

On every of \textquotedblleft random\textquotedblright\ steps
we take $d_{q}$ of randomly chosen uncovered points of
$PG(2,q)$ and compute the objective function $f$ adding each of
these $d_{q}$ points to $S_{i}$. The point providing the
extremum is included into $S_{i}.$ The value of $d_{q}$ is
given intuitively depending upon $q,$ upon the number of chosen
points (i.e. $|S_{i-1}|$), and upon the current task (small
arcs or the spectrum of arc sizes). For example, one can put
$d_{q}=1$ for finding of the spectrum and $ d_{q}=100\beta$
with $\beta=1,2,\ldots$ for small arcs.

As $S_{0}$ we can use a subset of points of an arc obtained in previous
stages of the search. Also, for finding the spectrum of arc sizes it is
fruitful to take as $S_{0}$ a part of points of a conic. A generator of
random numbers is used for a random choice. To get arcs with distinct sizes,
starting conditions of the generator are changed for the same set $S_{0}$.
In this way the algorithm works in a convenient limited region of the search
space to obtain examples improving the size of the arc from which the fixed
points have been taken.

In order to obtain arcs with new sizes one should make sufficiently many
attempts with the randomized greedy algorithms. For small arcs, the so
called predicted sizes considered in Section~\ref{sec4_observations} are
useful for understanding if a good result have been obtained. If the result
is not close to the predicted size, the attempts should be continued.

Note also that arcs with sizes close to $\overline{t}_{2}(2,q)$ usually are
obtained as a byproduct when we execute the computer search for the smallest
arcs using a few attempts.

\section{Small complete $k$-arcs in $PG(2,q)$, $q\leq 4561,$ $q\in T_{2}$}

\label{sec3_tables}

Throughout the paper, in all tables we denote
$A_{q}=\left\lfloor
a_{q}\sqrt{q}-\overline{t}_{2}(2,q)\right\rfloor $ where
\begin{eqnarray*}
a_{q}=\left\{
\begin{array}{cl}
4 & \text{if }q\leq 841,\text{ }q\in Q\\
4.5 & \text{if }853\leq q\leq 2593,\,q=2693,2753,\,q\notin Q \\
5 & \text{if }2609\le q\leq 8192, \,q\notin \{2693,2753\}.
\end{array}
\right. . \label{eq3_Aq_aq}
\end{eqnarray*}
Also, in all tables,
$B_{q}$ is a superior approximation of $\overline{t}_{2}(2,q)/
\sqrt{q}.$

For $q\leq 841$, the values of $\overline{t}_{2}(2,q)$ (up to
June 2009) are collected in \cite[Tab.\thinspace
1]{DFMP-JG2009}. In this work we obtained small arcs with new
sizes for $q\in T_{1}.$ The new arcs are obtained by computer
search, based on the randomized greedy algorithms. Complete
90-arcs for $q=601,661$ came from a theoretical study on orbits
of subgroups, helped by computer, see \cite{MaPa-private}.
These arcs are announced also in \cite[Tab.\thinspace
1]{DGMP-AMC}. A complete 104-arc for $q=709$ is obtained by the
greedy algorithm with the starting point set taking from
\cite{MaPa-private}. The current values of
$\overline{t}_{2}(2,q)$ for $q\leq 841$ are given in Table 1.
The data for $q\in T_{1}$ improving results of \cite[Tab.\thinspace
1]{DFMP-JG2009} are written in Table~1 in bold font. The exact
values $\overline{t} _{2}(2,q)=t_{2}(2,q)$ are marked by the
dot~\textquotedblleft$\centerdot$\textquotedblright. In
particular, due to the recent result \cite{MaMiP32} we noted
the values $t_{2}(2,31)=t_{2}(2,32)=14.$

\begin{figure}[tbp]
\textbf{Table 1}

The smallest known sizes $\overline{t}_{2}=\overline{t}_{2}(2,q)<4\sqrt{q}$
of complete arcs in planes $\mathrm{PG}(2,q),$ $q\leq 841.$

$A_{q}=\left\lfloor
4\sqrt{q}-\overline{t}_{2}(2,q)\right\rfloor $, $ B_{q}\geq
\overline{t}_{2}(2,q)/\sqrt{q}$\smallskip

\noindent$
\renewcommand{\arraystretch}{0.97}
\begin{array}{@{}r@{\,\,\,\,\,}l@{\,\,\,\,}r@{\,\,\,\,}c|c@{\,\,\,\,}c@{\,\,\,\,}r@{\,\,\,\,}c|
cc@{\,\,\,\,}cc|cr@{\,\,\,\,}rc@{}}
\hline
q & \overline{t}_{2} & A_{q} & B_{q} & q & \overline{t}_{2} & A_{q} & B_{q}
& q & \overline{t}_{2} & A_{q} & B_{q} & q & \overline{t}_{2} & A_{q} & B_{q}^{\phantom{H^{L}}}
\\ \hline
2 & \phantom{1}4\centerdot & 1 & 2.83 & 128 & 34 & 11 & 3.01 & 347 & 67 & 7
& 3.60 & \mathbf{599} & \mathbf{94} & \mathbf{3} & \mathbf{3.85} \\
3 & \phantom{1}4\centerdot & 2 & 2.31 & 131 & 36 & 9 & 3.15 & 349 & 67 & 7 &
3.59 & \mathbf{601} & \mathbf{90} & \mathbf{8} & \mathbf{3.68} \\
4 & \phantom{1}6\centerdot & 2 & 3.00 & 137 & 37 & 9 & 3.17 & 353 & 68 & 7 &
3.62 & \mathbf{607} & \mathbf{95} & \mathbf{3} & \mathbf{3.86} \\
5 & \phantom{1}6\centerdot & 2 & 2.69 & 139 & 37 & 10 & 3.14 & 359 & 69 & 6
& 3.65 & \mathbf{613} & \mathbf{96} & \mathbf{3} & \mathbf{3.88} \\
7 & \phantom{1}6\centerdot & 4 & 2.27 & 149 & 39 & 9 & 3.20 & 361 & 69 & 7 &
3.64 & \mathbf{617} & \mathbf{96} & \mathbf{3} & \mathbf{3.87} \\
8 & \phantom{1}6\centerdot & 5 & 2.13 & 151 & 39 & 10 & 3.18 & 367 & 70 & 6
& 3.66 & \mathbf{619} & \mathbf{96} & \mathbf{3} & \mathbf{3.86} \\
9 & \phantom{1}6\centerdot & 6 & 2.00 & 157 & 40 & 10 & 3.20 & \mathbf{373}
& \mathbf{70} & \mathbf{7} & \mathbf{3.63} & 625 & 96 & 4 & 3.84 \\
11 & \phantom{1}7\centerdot & 6 & 2.12 & 163 & 41 & 10 & 3.22 & 379 & 71 & 6
& 3.65 & \mathbf{631} & \mathbf{97} & \mathbf{3} & \mathbf{3.87} \\
13 & \phantom{1}8\centerdot & 6 & 2.22 & 167 & 42 & 9 & 3.26 & 383 & 71 & 7
& 3.63 & \mathbf{641} & \mathbf{98} & \mathbf{3} & \mathbf{3.88} \\
16 & \phantom{1}9\centerdot & 7 & 2.25 & 169 & 42 & 10 & 3.24 & 389 & 72 & 6
& 3.66 & 643 & 99 & 2 & 3.91 \\
17 & 10\centerdot & 6 & 2.43 & \mathbf{173} & \mathbf{43} & \mathbf{9} &
\mathbf{3.27} & 397 & 73 & 6 & 3.67 & 647 & 99 & 2 & 3.90 \\
19 & 10\centerdot & 7 & 2.30 & 179 & 44 & 9 & 3.29 & 401 & 74 & 6 & 3.70 &
653 & 100 & 2 & 3.92 \\
23 & 10\centerdot & 9 & 2.09 & \mathbf{181} & \mathbf{44} & \mathbf{9} &
\mathbf{3.28} & \mathbf{409} & \mathbf{74} & \mathbf{6} & \mathbf{3.66} & 659
& 100 & 2 & 3.90 \\
25 & 12\centerdot & 8 & 2.40 & 191 & 46 & 9 & 3.33 & 419 & 76 & 5 & 3.72 &
\mathbf{661} & \mathbf{90} & \mathbf{12} & \mathbf{3.51} \\
27 & 12\centerdot & 8 & 2.31 & \mathbf{193} & \mathbf{46} & \mathbf{9} &
\mathbf{3.32} & 421 & 76 & 6 & 3.71 & \mathbf{673} & \mathbf{101} & \mathbf{2
} & \mathbf{3.90} \\
29 & 13\centerdot & 8 & 2.42 & 197 & 47 & 9 & 3.35 & 431 & 77 & 6 & 3.71 &
\mathbf{677} & \mathbf{102} & \mathbf{2} & \mathbf{3.93} \\
31 & \mathbf{14\centerdot } & 8 & 2.52 & 199 & 47 & 9 & 3.34 & 433 & 77 & 6
& 3.71 & \mathbf{683} & \mathbf{102} & \mathbf{2} & \mathbf{3.91} \\
32 & \mathbf{14\centerdot } & 8 & 2.48 & 211 & 49 & 9 & 3.38 & 439 & 78 & 5
& 3.73 & \mathbf{691} & \mathbf{103} & \mathbf{2} & \mathbf{3.92} \\
37 & 15 & 9 & 2.47 & 223 & 51 & 8 & 3.42 & \mathbf{443} & \mathbf{78} &
\mathbf{6} & \mathbf{3.71} & 701 & 104 & 1 & 3.93 \\
41 & 16 & 9 & 2.50 & 227 & 51 & 9 & 3.39 & \mathbf{449} & \mathbf{79} &
\mathbf{5} & \mathbf{3.73} & \mathbf{709} & \mathbf{104} & \mathbf{2} & \mathbf{3.91} \\
43 & 16 & 10 & 2.45 & \mathbf{229} & \mathbf{51} & \mathbf{9} & \mathbf{3.38} & \mathbf{457} & \mathbf{80} &
\mathbf{5} & \mathbf{3.75} & 719 & 106 & 1 & 3.96 \\
47 & 18 & 9 & 2.63 & 233 & 52 & 9 & 3.41 & \mathbf{461} & \mathbf{80} &
\mathbf{5} & \mathbf{3.73} & 727 & 106 & 1 & 3.94 \\
49 & 18 & 10 & 2.58 & 239 & 53 & 8 & 3.43 & \mathbf{463} & \mathbf{80} &
\mathbf{6} & \mathbf{3.72} & 729 & 104 & 4 & 3.86 \\
53 & 18 & 11 & 2.48 & 241 & 53 & 9 & 3.42 & \mathbf{467} & \mathbf{81} &
\mathbf{5} & \mathbf{3.75} & 733 & 107 & 1 & 3.96 \\
59 & 20 & 10 & 2.61 &  \mathbf{243} &  \mathbf{53} &  \mathbf{9} &  \mathbf{3.4}0 & \mathbf{479} & \mathbf{82} &
\mathbf{5} & \mathbf{3.75} & 739 & 107 & 1 & 3.94 \\
61 & 20 & 11 & 2.57 & 251 & 55 & 8 & 3.48 & \mathbf{487} & \mathbf{83} &
\mathbf{5} & \mathbf{3.77} & 743 & 108 & 1 & 3.97 \\
64 & 22 & 10 & 2.75 & 256 & 55 & 9 & 3.44 & \mathbf{491} & \mathbf{83} &
\mathbf{5} & \mathbf{3.75} & 751 & 108 & 1 & 3.95 \\
67 & 23 & 9 & 2.81 & \mathbf{257} & \mathbf{55} & \mathbf{9} & \mathbf{3.44}
& \mathbf{499} & \mathbf{84} & \mathbf{5} & \mathbf{3.77} & 757 & 109 & 1 &
3.97 \\
71 & 22 & 11 & 2.62 & 263 & 56 & 8 & 3.46 & 503 & 85 & 4 & 3.79 & 761 & 109
& 1 & 3.96 \\
73 & 24 & 10 & 2.81 & 269 & 57 & 8 & 3.48 & 509 & 85 & 5 & 3.77 & 769 & 110
& 0 & 3.97 \\
79 & 26 & 9 & 2.93 & \mathbf{271} & \mathbf{57} & \mathbf{8} & \mathbf{3.47}
& 512 & 86 & 4 & 3.81 & 773 & 111 & 0 & 4.00 \\
81 & 26 & 10 & 2.89 & \mathbf{277} & \mathbf{58} & \mathbf{8} & \mathbf{3.49}
& 521 & 86 & 5 & 3.77 & 787 & 112 & 0 & 4.00 \\
83 & 27 & 9 & 2.97 & 281 & 59 & 8 & 3.52 & 523 & 86 & 5 & 3.77 & 797 & 112 &
0 & 3.97 \\
89 & 28 & 9 & 2.97 & 283 & 59 & 8 & 3.51 & \mathbf{529} & \mathbf{87} & \mathbf{5} & \mathbf{3.79} & 809 & 113 &
0 & 3.98 \\
97 & 30 & 9 & 3.05 & 289 & 60 & 8 & 3.53 & 541 & 89 & 4 & 3.83 & 811 & 113 &
0 & 3.97 \\
101 & 30 & 10 & 2.99 & \mathbf{293} & \mathbf{60} & \mathbf{8} & \mathbf{3.51
} & 547 & 89 & 4 & 3.81 & 821 & 114 & 0 & 3.98 \\
103 & 31 & 9 & 3.06 & 307 & 62 & 8 & 3.54 & 557 & 90 & 4 & 3.82 & 823 & 114
& 0 & 3.98 \\
107 & 32 & 9 & 3.10 & 311 & 63 & 7 & 3.58 & \mathbf{563} & \mathbf{91} &
\mathbf{3} & \mathbf{3.84} & 827 & 115 & 0 & 4.00 \\
109 & 32 & 9 & 3.07 & 313 & 63 & 7 & 3.57 & \mathbf{569} & \mathbf{91} &
\mathbf{4} & \mathbf{3.82} & 829 & 115 & 0 & 4.00 \\
113 & 33 & 9 & 3.11 & 317 & 63 & 8 & 3.54 & \mathbf{571} & \mathbf{92} &
\mathbf{3} & \mathbf{3.86} & 839 & 115 & 0 & 3.98 \\
121 & 34 & 10 & 3.10 & 331 & 65 & 7 & 3.58 & \mathbf{577} & \mathbf{92} &
\mathbf{4} & \mathbf{3.84} & 841 & 112 & 4 & 3.87 \\
125 & 35 & 9 & 3.14 & 337 & 66 & 7 & 3.60 & \mathbf{587} & \mathbf{93} &
\mathbf{3} & \mathbf{3.84} &  &  &  &  \\
127 & 35 & 10 & 3.11 & \mathbf{343} & \mathbf{66} & \mathbf{8} & \mathbf{3.57
} & \mathbf{593} & \mathbf{94} & \mathbf{3} & \mathbf{3.87} &  &  &  &  \\
\hline
\end{array}
$
\end{figure}

From Table 1 and the results of \cite{Giul2000,GiulUghi}, on
complete $(4\sqrt{q}-4)$-arcs for $q=p^{2}$ (see Introduction)
we obtain Theorem \ref {Th3_4sqroot(q)} improving and extending
the results of \cite[Th.\thinspace 1]{DFMP-JG2009}.

\begin{thm}
\label{Th3_4sqroot(q)} In $PG(2,q),$ the following holds.
\begin{eqnarray}
t_{2}(2,q)&<&\phantom{.5}4\sqrt{q}~ \mbox{ for }2\leq q\leq 841,\text{ }
q\in Q.  \label{eq3_<4sqroot(q)}\\
t_{2}(2,q) &\leq &\phantom{.5}3\sqrt{q}~\mbox{ for }2\leq q\leq 89,\text{ }
\,~q=101;  \notag \\
t_{2}(2,q) &<&3.5\sqrt{q}~\mbox{ for }2\leq q\leq 277;
\label{eq3_4sqroot(q)_a*sqroot(q)} \notag\\
t_{2}(2,q) &<&3.8\sqrt{q}~\mbox{ for }2\leq q\leq 509,~q=521,523,529,601,661.
\notag
\end{eqnarray}
 Also,
\begin{eqnarray*}
t_{2}(2,q) &\leq &4\sqrt{q}-9\ \mbox{ for }\,\,37\leq q\leq 211,\text{ }
q=23,227,229,233,241,243,256,257,661;  \notag \\
t_{2}(2,q) &\leq &4\sqrt{q}-8\ \mbox{ for }\,\,23\leq q\leq 307,\text{ }
q=317,343,601,661;  \notag\\
t_{2}(2,q) &\leq &4\sqrt{q}-7\ \mbox{ for }\ 19\leq q\leq 353,\
q=16,361,373,383,601,661;  \notag\\
t_{2}(2,q) &\leq &4\sqrt{q}-6\ \mbox{ for
}\,\,\phantom{1}9\leq q\leq 409,\ q=421,431,433,443,463,601,661;
\label{eq3_4sqroot(q)-delta}
\end{eqnarray*}
\begin{eqnarray*}
t_{2}(2,q) &\leq &4\sqrt{q}-5\ \mbox{ for }\ \phantom{1}8\leq q\leq 499,\
q=509,521,523,529,601,661;  \notag \\
t_{2}(2,q) &\leq &4\sqrt{q}-4\ \mbox{ for }\ \phantom{1}7\leq q\leq 557,\
q=569,577,601,625,661,729,841,\,q\in Q;\notag\\
t_{2}(2,q) &< &4\sqrt{q}-3\ \mbox{ for }\ \phantom{1}7\leq q\leq 641,\
q=661,729,841,\,q\in Q;  \notag \\
t_{2}(2,q) &\leq &4\sqrt{q}-2\ \mbox{ for }\ \phantom{1}3\leq q\leq 691,
\text{ }q=709,729,841,\,q\in Q;  \notag \\
t_{2}(2,q) &< &4\sqrt{q}-1\ \mbox{ for }\ \phantom{1}2\leq q\leq 761,
\text{ }q=841,\,q\in Q.  \notag
\end{eqnarray*}
\end{thm}

In Table 2, the current values of $\overline{t}_{2}(2,q)$ for
$853\leq q\leq 2593$ are given. The data for $q=p^{2}$ with
$\overline{t}_{2}(2,q)=4\sqrt{q} -4$ \cite{Giul2000,GiulUghi}
are written in bold font.
\begin{figure}[tbp]
\noindent \textbf{Table 2}

\noindent The smallest known sizes $\overline{t}_{2}=\overline{t}
_{2}(2,q)<4.5\sqrt{q}$ of complete arcs in planes $PG(2,q),$

\noindent $853\leq q\leq 2593,$ $A_{q}=\left\lfloor
a_{q}\sqrt{q}-\overline{t }_{2}(2,q)\right\rfloor $, $B_{q}\geq
\overline{t}_{2}(2,q)/\sqrt{q}$\smallskip

\noindent$\renewcommand{\arraystretch}{1.0}
\begin{array}{@{}r@{\,\,\,}c@{\,\,\,}c@{\,\,\,\,}c@{\,\,}|@{\,\,}c@{\,\,\,}
c@{\,\,\,\,}c@{\,\,\,}c@{\,\,}|@{\,\,}c@{\,\,\,\,}c@{\,\,\,\,}c@{\,\,\,}c@{\,\,}|@{\,\,}c@{\,\,\,\,}
c@{\,\,\,\,}c@{\,\,\,\,}c@{}}
\hline
q & \overline{t}_{2} & A_{q} & B_{q} & q & \overline{t}_{2} &
A_{q} & B_{q} & q & \overline{t}_{2} & A_{q} & B_{q} & q &
\overline{t}_{2} & A_{q} & B_{q}^{\phantom{H^{L}}} \\ \hline
853 & 117 & 14 & 4.01 & 1277 & 150 & 10 & 4.20 & 1693 & 178 & 7
& 4.33 & 2141
& 205 & 3 & 4.44 \\
857 & 118 & 13 & 4.04 & 1279 & 150 & 10 & 4.20 & 1697 & 178 & 7
& 4.33 & 2143
& 205 & 3 & 4.43 \\
859 & 118 & 13 & 4.03 & 1283 & 150 & 11 & 4.19 & 1699 & 178 & 7
& 4.32 & 2153
& 205 & 3 & 4.42 \\
863 & 118 & 14 & 4.02 & 1289 & 151 & 10 & 4.21 & 1709 & 178 & 8
& 4.31 & 2161
& 205 & 4 & 4.41 \\
877 & 119 & 14 & 4.02 & 1291 & 151 & 10 & 4.21 & 1721 & 180 & 6
& 4.34 & 2179
& 207 & 3 & 4.44 \\
881 & 120 & 13 & 4.05 & 1297 & 151 & 11 & 4.20 & 1723 & 180 & 6
& 4.34 & 2187
& 207 & 3 & 4.43 \\
883 & 120 & 13 & 4.04 & 1301 & 152 & 10 & 4.22 & 1733 & 180 & 7
& 4.33 & 2197
& 208 & 2 & 4.44 \\
887 & 120 & 14 & 4.03 & 1303 & 151 & 11 & 4.19 & 1741 & 181 & 6
& 4.34 & 2203
& 208 & 3 & 4.44 \\
907 & 122 & 13 & 4.06 & 1307 & 152 & 10 & 4.21 & 1747 & 181 & 7
& 4.34 & 2207
& 208 & 3 & 4.43 \\
911 & 122 & 13 & 4.05 & 1319 & 153 & 10 & 4.22 & 1753 & 182 & 6
& 4.35 & 2209
& 208 & 3 & 4.43 \\
919 & 123 & 13 & 4.06 & 1321 & 153 & 10 & 4.21 & 1759 & 182 & 6
& 4.34 & 2213
& 209 & 2 & 4.45 \\
929 & 124 & 13 & 4.07 & 1327 & 153 & 10 & 4.21 & 1777 & 183 & 6
& 4.35 & 2221
& 209 & 3 & 4.44 \\
937 & 124 & 13 & 4.06 & 1331 & 154 & 10 & 4.23 & 1783 & 183 & 7
& 4.34 & 2237
& 210 & 2 & 4.45 \\
941 & 125 & 13 & 4.08 & 1361 & 156 & 10 & 4.23 & 1787 & 183 & 7
& 4.33 & 2239
& 210 & 2 & 4.44 \\
947 & 125 & 13 & 4.07 & 1367 & 156 & 10 & 4.22 & 1789 & 184 & 6
& 4.36 & 2243
& 210 & 3 & 4.44 \\
953 & 126 & 12 & 4.09 & \mathbf{1369} & \mathbf{144} & \mathbf{4} & \mathbf{3.90}
 & 1801 & 184 & 6 & 4.34 & 2251 & 211 & 2 & 4.45 \\
\mathbf{961} & \mathbf{120} & \mathbf{4} & \mathbf{3.88} & 1373
& 157 & 9 &
4.24 & 1811 & 184 & 7 & 4.33 & 2267 & 211 & 3 & 4.44 \\
967 & 127 & 12 & 4.09 & 1381 & 157 & 10 & 4.23 & 1823 & 186 & 6
& 4.36 & 2269
& 212 & 2 & 4.46 \\
971 & 127 & 13 & 4.08 & 1399 & 159 & 9 & 4.26 & 1831 & 186 & 6
& 4.35 & 2273
& 212 & 2 & 4.45 \\
977 & 127 & 13 & 4.07 & 1409 & 159 & 9 & 4.24 & 1847 & 187 & 6
& 4.36 & 2281
& 213 & 1 & 4.46 \\
983 & 128 & 13 & 4.09 & 1423 & 160 & 9 & 4.25 & 1849 & 187 & 6
& 4.35 & 2287
& 213 & 2 & 4.46 \\
991 & 127 & 14 & 4.04 & 1427 & 160 & 9 & 4.24 & 1861 & 188 & 6
& 4.36 & 2293
& 213 & 2 & 4.45 \\
997 & 129 & 13 & 4.09 & 1429 & 161 & 9 & 4.26 & 1867 & 189 & 5
& 4.38 & 2297
& 213 & 2 & 4.45 \\
1009 & 130 & 12 & 4.10 & 1433 & 161 & 9 & 4.26 & 1871 & 189 & 5
& 4.37 & 2309
& 214 & 2 & 4.46 \\
1013 & 130 & 13 & 4.09 & 1439 & 161 & 9 & 4.25 & 1873 & 189 & 5
& 4.37 & 2311
& 214 & 2 & 4.46 \\
1019 & 131 & 12 & 4.11 & 1447 & 162 & 9 & 4.26 & 1877 & 189 & 5
& 4.37 & 2333
& 215 & 2 & 4.46 \\
1021 & 131 & 12 & 4.10 & 1451 & 162 & 9 & 4.26 & 1879 & 189 & 6
& 4.37 & 2339
& 216 & 1 & 4.47 \\
\mathbf{1024} & \mathbf{124} & \mathbf{4} & \mathbf{3.88} &
1453 & 162 & 9 &
4.25 & 1889 & 190 & 5 & 4.38 & 2341 & 216 & 1 & 4.47 \\
1031 & 132 & 12 & 4.12 & 1459 & 163 & 8 & 4.27 & 1901 & 191 & 5
& 4.39 & 2347
& 216 & 2 & 4.46 \\
1033 & 132 & 12 & 4.11 & 1471 & 163 & 9 & 4.25 & 1907 & 191 & 5
& 4.38 & 2351
& 216 & 2 & 4.46 \\
1039 & 132 & 13 & 4.10 & 1481 & 164 & 9 & 4.27 & 1913 & 191 & 5
& 4.37 & 2357
& 217 & 1 & 4.47 \\
1049 & 133 & 12 & 4.11 & 1483 & 164 & 9 & 4.26 & 1931 & 192 & 5
& 4.37 & 2371
& 217 & 2 & 4.46 \\
1051 & 133 & 12 & 4.11 & 1487 & 164 & 9 & 4.26 & 1933 & 192 & 5
& 4.37 & 2377
& 216 & 3 & 4.44 \\
1061 & 134 & 12 & 4.12 & 1489 & 165 & 8 & 4.28 & 1949 & 193 & 5
& 4.38 & 2381
& 217 & 2 & 4.45 \\
1063 & 134 & 12 & 4.11 & 1493 & 165 & 8 & 4.28 & 1951 & 194 & 4
& 4.40 & 2383
& 218 & 1 & 4.47 \\
1069 & 135 & 12 & 4.13 & 1499 & 165 & 9 & 4.27 & 1973 & 195 & 4
& 4.40 & 2389
& 218 & 1 & 4.47 \\
1087 & 136 & 12 & 4.13 & 1511 & 166 & 8 & 4.28 & 1979 & 195 & 5
& 4.39 & 2393
& 218 & 2 & 4.46 \\
1091 & 136 & 12 & 4.12 & 1523 & 167 & 8 & 4.28 & 1987 & 196 & 4
& 4.40 & 2399
& 219 & 1 & 4.48 \\
1093 & 136 & 12 & 4.12 & 1531 & 167 & 9 & 4.27 & 1993 & 196 & 4
& 4.40 &
\mathbf{2401} & \mathbf{192} & \mathbf{4} & \mathbf{3.92} \\
1097 & 137 & 12 & 4.14 & 1543 & 167 & 9 & 4.26 & 1997 & 196 & 5
& 4.39 & 2411
& 220 & 0 & 4.49 \\
1103 & 137 & 12 & 4.13 & 1549 & 169 & 8 & 4.30 & 1999 & 196 & 5
& 4.39 & 2417
& 220 & 1 & 4.48 \\
1109 & 138 & 11 & 4.15 & 1553 & 169 & 8 & 4.29 & 2003 & 197 & 4
& 4.41 & 2423 & 220 & 1 & 4.47 \\ \hline
\end{array}
$
\end{figure}

\begin{figure}[t]
\noindent \textbf{Table 2 }(continue)

\noindent The smallest known sizes $\overline{t}_{2}=\overline{t}
_{2}(2,q)<4.5\sqrt{q}$ of complete arcs in planes $PG(2,q),$

\noindent $853\leq q\leq 2593,$ $A_{q}=\left\lfloor
a_{q}\sqrt{q}-\overline{t} _{2}(2,q)\right\rfloor $, $B_{q}\geq
\overline{t}_{2}(2,q)/\sqrt{q}$\medskip

\noindent$\renewcommand{\arraystretch}{1.0}
\begin{array}{@{}r@{\,\,\,\,}c@{\,\,\,\,}c@{\,\,\,\,}c|@{\,\,}
c@{\,\,\,\,}c@{\,\,\,\,}c@{\,\,\,}c|@{\,\,}c@{\,\,\,\,}c@{\,\,\,\,}c@{\,\,\,}c|
@{\,\,}c@{\,\,\,\,}c@{\,\,\,\,}c@{\,\,\,\,}c@{}}
\hline
q & \overline{t}_{2} & A_{q} & B_{q} & q & \overline{t}_{2} & A_{q} & B_{q}
& q & \overline{t}_{2} & A_{q} & B_{q} & q & \overline{t}_{2} & A_{q} &
B_{q}^{\phantom{H^{L}}} \\ \hline
1117 & 138 & 12 & 4.13 & 1559 & 169 & 8 & 4.29 & 2011 & 197 & 4 & 4.40 & 2437
& 221 & 1 & 4.48 \\
1123 & 139 & 11 & 4.15 & 1567 & 170 & 8 & 4.30 & 2017 & 197 & 5 & 4.39 & 2441
& 221 & 1 & 4.48 \\
1129 & 139 & 12 & 4.14 & 1571 & 170 & 8 & 4.29 & 2027 & 198 & 4 & 4.40 & 2447
& 221 & 1 & 4.47 \\
1151 & 141 & 11 & 4.16 & 1579 & 170 & 8 & 4.28 & 2029 & 198 & 4 & 4.40 & 2459
& 222 & 1 & 4.48 \\
1153 & 141 & 11 & 4.16 & 1583 & 171 & 8 & 4.30 & 2039 & 199 & 4 & 4.41 & 2467
& 223 & 0 & 4.49 \\
1163 & 142 & 11 & 4.17 & 1597 & 172 & 7 & 4.31 & 2048 & 199 & 4 & 4.40 & 2473
& 223 & 0 & 4.49 \\
1171 & 142 & 11 & 4.15 & 1601 & 172 & 8 & 4.30 & 2053 & 200 & 3 & 4.42 & 2477
& 223 & 0 & 4.49 \\
1181 & 143 & 11 & 4.17 & 1607 & 172 & 8 & 4.30 & 2063 & 200 & 4 & 4.41 & 2503
& 225 & 0 & 4.50 \\
1187 & 144 & 11 & 4.18 & 1609 & 172 & 8 & 4.29 & 2069 & 200 & 4 & 4.40 & 2521
& 225 & 0 & 4.49 \\
1193 & 144 & 11 & 4.17 & 1613 & 173 & 7 & 4.31 & 2081 & 201 & 4 & 4.41 & 2531
& 226 & 0 & 4.50 \\
1201 & 145 & 10 & 4.19 & 1619 & 173 & 8 & 4.30 & 2083 & 201 & 4 & 4.41 & 2539
& 226 & 0 & 4.49 \\
1213 & 145 & 11 & 4.17 & 1621 & 173 & 8 & 4.30 & 2087 & 201 & 4 & 4.40 & 2543
& 226 & 0 & 4.49 \\
1217 & 146 & 10 & 4.19 & 1627 & 174 & 7 & 4.32 & 2089 & 202 & 3 & 4.42 & 2549
& 226 & 1 & 4.48 \\
1223 & 146 & 11 & 4.18 & 1637 & 174 & 8 & 4.31 & 2099 & 202 & 4 & 4.41 & 2551
& 227 & 0 & 4.50 \\
1229 & 146 & 11 & 4.17 & 1657 & 175 & 8 & 4.30 & 2111 & 203 & 3 & 4.42 & 2557
& 227 & 0 & 4.49 \\
1231 & 147 & 10 & 4.19 & 1663 & 176 & 7 & 4.32 & 2113 & 203 & 3 & 4.42 & 2579
& 228 & 0 & 4.49 \\
1237 & 147 & 11 & 4.18 & 1667 & 176 & 7 & 4.32 & 2129 & 204 & 3 & 4.43 & 2591
& 229 & 0 & 4.50 \\
1249 & 148 & 11 & 4.19 & 1669 & 176 & 7 & 4.31 & 2131 & 204 & 3 & 4.42 & 2593
& 229 & 0 & 4.50 \\
1259 & 149 & 10 & 4.20 & \mathbf{1681} & \mathbf{160} & \mathbf{4} & \mathbf{3.91}
& 2137 & 204 & 4 & 4.42 &  &  &  &  \\ \hline
\end{array}
$
\end{figure}

From Table 2, we obtain Theorem \ref{Th3_4.5sqroot(q)}.

\begin{thm}
\label{Th3_4.5sqroot(q)} In $PG(2,q),$ the following  holds.
\begin{eqnarray}
t_{2}(2,q)&<&4.5\sqrt{q}~ \mbox{ for }q\leq 2593,\,q=2693,2753.
\label{eq3_<4.5sqroot(q)}\\
t_{2}(2,q) &<&4.1\sqrt{q}~\mbox{ for }q\leq
1013,\,q=1021,1024,1039,1369,1681,2401;  \notag \\
t_{2}(2,q) &<&4.2\sqrt{q}~\mbox{ for }q\leq 1283,\,q=1297,1303,1369,1681,2401;
\notag \\
t_{2}(2,q) &<&4.3\sqrt{q}~\mbox{ for }q\leq
1583,\,q=1601,1607,1609,1619,1621,1657,1681,2401;
\label{eq3_4.5sqroot(q)_a*sqroot(q)} \notag\\
t_{2}(2,q) &<&4.4\sqrt{q}~\mbox{ for
}q\leq 1999,\,q=2011,2017,2027,2029,2048,2069,2087,2401.  \notag
\end{eqnarray}
Also,
\begin{eqnarray*}
t_{2}(2,q) &<&4.5\sqrt{q}-12\,\mbox{ for }q\leq 1103,\text{\thinspace }
q=1117,1129,1369,1681,2401;  \notag \\
t_{2}(2,q) &<&4.5\sqrt{q}-11\,\mbox{ for }q\leq
1193,\,q=1213,1223,1229,1237,1249,1283,1297,1303,\\
&&\phantom{4.5\sqrt{q}-11\,\mbox{ for }q\leq
1193,\,q=}\,\,1369,1681,2401;\notag\\
t_{2}(2,q) &<&4.5\sqrt{q}-10\,\mbox{ for }q\leq
1369,\,q=1381,1681,2401;  \notag \\
t_{2}(2,q) &<&4.5\sqrt{q}-9\,~\,\mbox{ for }q\leq 1453,\,q=
1471,1481,1483,1487,1499,1531,1543,1681,\notag\\
&&\phantom{4.5\sqrt{q}-9\,~\,\mbox{ for }q\leq 1453,\,q=}~2401;\notag \\
t_{2}(2,q) &<&4.5\sqrt{q}-8\,~\,\mbox{ for }q\leq
1583,\,q=1601,1607,1609,1619,1621,1637,1657,1681,\notag \\
&&\phantom{4.5\sqrt{q}-8\,~\,\mbox{ for }q\leq
1571,\,q=}~1709,2401; \notag\\
t_{2}(2,q) &<&4.5\sqrt{q}-7\,~\,\mbox{ for }q\leq
1709,\,q=1733,1747,1783,1787,1811,2401;  \label{eq3_4.5sqroot(q)--delta} \\
t_{2}(2,q) &<&4.5\sqrt{q}-6\,~\,\mbox{ for
}q\leq 1861,\,q=1879,2401;  \notag\\
t_{2}(2,q) &<&4.5\sqrt{q}-5\,~\,\mbox{ for
}q\leq 1949,\,q=1979,1997,1999,2017,2401;  \notag\\
t_{2}(2,q) &<&4.5\sqrt{q}-4\,~\,\mbox{ for
}q\leq 2048,\,q=2063,2069,2081,2083,2087,2099,2137,2161,\notag \\
&&\phantom{4.5\sqrt{q}-8\,~\,\mbox{ for }q\leq
1571,\,q=}~2401;  \notag \\
t_{2}(2,q) &<&4.5\sqrt{q}-3\,~\,\mbox{ for }q\leq
2187,\,q=2203,2207,2209,2221,2243,2267,2377,2401;  \notag \\
t_{2}(2,q) &<&4.5\sqrt{q}-2\,~\,\mbox{ for }q\leq
2273,\,q=2287,2293,2297,2309,2311,2333,2347,2351,\notag \\
&&\phantom{4.5\sqrt{q}-2\,~\,\mbox{ for }q\leq
2273,\,q=\,\,}2371,2377,2381,2393,2401;  \notag \\
t_{2}(2,q) &<&4.5\sqrt{q}-1\,~\mbox{ for }q\leq
2401,\,q=2417,2423,2437,2441,2447,2459,2549.  \notag
\end{eqnarray*}
\end{thm}

In Table 3, the current values of $\overline{t}_{2}(2,q)$  for
$2609\leq q\leq 4561$  are given.  The data with
$\overline{t}_{2}(2,q)<4.5\sqrt{q}$ are written in bold font.

Values of $\overline{t}_{2}(2,q)$ for relatively great $q\in
T_{2}\cup T_{3}$ are given in Table 4. The notation
$\overline{D}_{q}(\frac{3}{4})$ is explained in the next
section.

In Table 2 for $ q\in Q,$ we use the results of
\cite{DGMP-JCD,Giul2000,GiulUghi}, see also \cite[p.\thinspace
35]{DFMP-JG2009}. In Table 4, for $q\in T_{3}$ we use the
results of \cite {DGMP-Innov,DGMP-JCD}, see also
\cite[p.\thinspace 35]{DFMP-JG2009} and Introduction. The rest
of sizes $k$ for small complete $k$-arcs in Tables 2, 3 and 4
is obtained in this work by computer search with the help of
the randomized greedy algorithms.

Note that a complete $199$-arc in $PG(2,2048)$ of Table 2, a
complete $301$-arc in $PG(2,4096)$ of Table 3, and a complete
$450$-arc in $PG(2,8192)$ of Table 4 improve the results of
\cite {DGMP-JCD} for $q=2^{11},2^{12},2^{13}$, see
Introduction.

\begin{figure}[tbp]
\noindent \textbf{Table 3}

\noindent The smallest known sizes
$\overline{t}_{2}=\overline{t} _{2}(2,q)<4.75\sqrt{q}$ of
complete arcs in planes $PG(2,q),$

\noindent $2609\leq q\leq 4561$, $A_{q}=\left\lfloor
a_{q}\sqrt{q}-\overline{ t}_{2}(2,q)\right\rfloor $, $B_{q}\geq
\overline{t}_{2}(2,q)/\sqrt{q}$\smallskip

\noindent$\renewcommand{\arraystretch}{1.0}
\begin{array}{@{}r@{\,\,\,\,}c@{\,\,\,\,}c@{\,\,\,}c|@{\,\,}c@{\,\,\,\,}
c@{\,\,\,\,}c@{\,\,\,}c|@{\,\,}c@{\,\,\,\,}c@{\,\,\,\,}c@{\,\,\,}c|
@{\,\,}c@{\,\,\,\,}c@{\,\,\,\,}c@{\,\,\,\,}c@{}}
\hline
q & \overline{t}_{2} & A_{q} & B_{q} & q & \overline{t}_{2} & A_{q} & B_{q}
& q & \overline{t}_{2} & A_{q} & B_{q} & q & \overline{t}_{2} & A_{q} &
B_{q}^{\phantom{H^{L}}} \\ \hline
2609 & 230 & 25 & 4.51 & 3079 & 253 & 24 & 4.56 & 3571 & 277 & 21 & 4.64 &
4073 & 299 & 20 & 4.69 \\
2617 & 231 & 24 & 4.52 & 3083 & 253 & 24 & 4.56 & 3581 & 278 & 21 & 4.65 &
4079 & 300 & 19 & 4.70 \\
2621 & 231 & 24 & 4.52 & 3089 & 254 & 23 & 4.58 & 3583 & 277 & 22 & 4.63 &
4091 & 300 & 19 & 4.70 \\
2633 & 231 & 25 & 4.51 & 3109 & 255 & 23 & 4.58 & 3593 & 278 & 21 & 4.64 &
4093 & 300 & 19 & 4.69 \\
2647 & 232 & 25 & 4.51 & 3119 & 255 & 24 & 4.57 & 3607 & 278 & 22 & 4.63 &
4096 & 301 & 19 & 4.71 \\
2657 & 233 & 24 & 4.53 & 3121 & 255 & 24 & 4.57 & 3613 & 279 & 21 & 4.65 &
4099 & 300 & 20 & 4.69 \\
2659 & 233 & 24 & 4.52 & 3125 & 256 & 23 & 4.58 & 3617 & 278 & 22 & 4.63 &
4111 & 301 & 19 & 4.70 \\
2663 & 233 & 25 & 4.52 & 3137 & 257 & 23 & 4.59 & 3623 & 279 & 21 & 4.64 &
4127 & 302 & 19 & 4.71 \\
2671 & 233 & 25 & 4.51 & 3163 & 257 & 24 & 4.57 & 3631 & 280 & 21 & 4.65 &
4129 & 302 & 19 & 4.70 \\
2677 & 234 & 24 & 4.53 & 3167 & 258 & 23 & 4.59 & 3637 & 280 & 21 & 4.65 &
4133 & 302 & 19 & 4.70 \\
2683 & 234 & 24 & 4.52 & 3169 & 258 & 23 & 4.59 & 3643 & 281 & 20 & 4.66 &
4139 & 302 & 19 & 4.70 \\
2687 & 234 & 25 & 4.52 & 3181 & 259 & 23 & 4.60 & 3659 & 281 & 21 & 4.65 &
4153 & 301 & 21 & 4.68 \\
2689 & 234 & 25 & 4.52 & 3187 & 258 & 24 & 4.58 & 3671 & 282 & 20 & 4.66 &
4157 & 303 & 19 & 4.70 \\
\mathbf{2693} & \mathbf{233} & \mathbf{0} & \mathbf{4.49} & 3191 & 259 & 23
& 4.59 & 3673 & 282 & 21 & 4.66 & 4159 & 302 & 20 & 4.69 \\
2699 & 235 & 24 & 4.53 & 3203 & 259 & 23 & 4.58 & 3677 & 282 & 21 & 4.66 &
4177 & 303 & 20 & 4.69 \\
2707 & 235 & 25 & 4.52 & 3209 & 260 & 23 & 4.59 & 3691 & 283 & 20 & 4.66 &
4201 & 305 & 19 & 4.71 \\
2711 & 235 & 25 & 4.52 & 3217 & 261 & 22 & 4.61 & 3697 & 283 & 21 & 4.66 &
4211 & 305 & 19 & 4.71 \\
2713 & 236 & 24 & 4.54 & 3221 & 261 & 22 & 4.60 & 3701 & 283 & 21 & 4.66 &
4217 & 305 & 19 & 4.70 \\
2719 & 236 & 24 & 4.53 & 3229 & 260 & 24 & 4.58 & 3709 & 284 & 20 & 4.67 &
4219 & 305 & 19 & 4.70 \\
2729 & 236 & 25 & 4.52 & 3251 & 262 & 23 & 4.60 & 3719 & 284 & 20 & 4.66 &
4229 & 305 & 20 & 4.70 \\
2731 & 237 & 24 & 4.54 & 3253 & 261 & 24 & 4.58 & 3721 & 284 & 21 & 4.66 &
4231 & 306 & 19 & 4.71 \\
2741 & 237 & 24 & 4.53 & 3257 & 263 & 22 & 4.61 & 3727 & 284 & 21 & 4.66 &
4241 & 306 & 19 & 4.70 \\
2749 & 237 & 25 & 4.53 & 3259 & 263 & 22 & 4.61 & 3733 & 284 & 21 & 4.65 &
4243 & 306 & 19 & 4.70 \\
\mathbf{2753} & \mathbf{236} & \mathbf{0} & \mathbf{4.50} & 3271 & 263 & 22
& 4.60 & 3739 & 283 & 22 & 4.63 & 4253 & 306 & 20 & 4.70 \\
2767 & 238 & 25 & 4.53 & 3299 & 264 & 23 & 4.60 & 3761 & 286 & 20 & 4.67 &
4259 & 307 & 19 & 4.71 \\
2777 & 238 & 25 & 4.52 & 3301 & 264 & 23 & 4.60 & 3767 & 286 & 20 & 4.66 &
4261 & 307 & 19 & 4.71 \\
2789 & 240 & 24 & 4.55 & 3307 & 265 & 22 & 4.61 & 3769 & 286 & 20 & 4.66 &
4271 & 307 & 19 & 4.70 \\
2791 & 240 & 24 & 4.55 & 3313 & 265 & 22 & 4.61 & 3779 & 286 & 21 & 4.66 &
4273 & 306 & 20 & 4.69 \\
2797 & 240 & 24 & 4.54 & 3319 & 266 & 22 & 4.62 & 3793 & 287 & 20 & 4.67 &
4283 & 308 & 19 & 4.71 \\
2801 & 240 & 24 & 4.54 & 3323 & 265 & 23 & 4.60 & 3797 & 287 & 21 & 4.66 &
4289 & 308 & 19 & 4.71 \\
2803 & 240 & 24 & 4.54 & 3329 & 266 & 22 & 4.62 & 3803 & 287 & 21 & 4.66 &
4297 & 308 & 19 & 4.70 \\
2809 & 241 & 24 & 4.55 & 3331 & 266 & 22 & 4.61 & 3821 & 288 & 21 & 4.66 &
4327 & 310 & 18 & 4.72 \\
2819 & 241 & 24 & 4.54 & 3343 & 265 & 24 & 4.59 & 3823 & 289 & 20 & 4.68 &
4337 & 311 & 18 & 4.73 \\
2833 & 242 & 24 & 4.55 & 3347 & 266 & 23 & 4.60 & 3833 & 289 & 20 & 4.67 &
4339 & 311 & 18 & 4.73 \\
2837 & 242 & 24 & 4.55 & 3359 & 267 & 22 & 4.61 & 3847 & 290 & 20 & 4.68 &
4349 & 311 & 18 & 4.72 \\
2843 & 242 & 24 & 4.54 & 3361 & 267 & 22 & 4.61 & 3851 & 290 & 20 & 4.68 &
4357 & 311 & 19 & 4.72 \\
2851 & 243 & 23 & 4.56 & 3371 & 268 & 22 & 4.62 & 3853 & 289 & 21 & 4.66 &
4363 & 310 & 20 & 4.70 \\
2857 & 243 & 24 & 4.55 & 3373 & 268 & 22 & 4.62 & 3863 & 290 & 20 & 4.67 &
4373 & 312 & 18 & 4.72 \\
2861 & 243 & 24 & 4.55 & 3389 & 268 & 23 & 4.61 & 3877 & 291 & 20 & 4.68 &
4391 & 313 & 18 & 4.73 \\
2879 & 244 & 24 & 4.55 & 3391 & 269 & 22 & 4.62 & 3881 & 291 & 20 & 4.68 &
4397 & 313 & 18 & 4.73 \\
2887 & 243 & 25 & 4.53 & 3407 & 270 & 21 & 4.63 & 3889 & 291 & 20 & 4.67 &
4409 & 314 & 18 & 4.73 \\
2897 & 245 & 24 & 4.56 & 3413 & 269 & 23 & 4.61 & 3907 & 292 & 20 & 4.68 &
4421 & 314 & 18 & 4.73 \\ \hline
\end{array}
$
\end{figure}

\begin{figure}[tb]
\noindent \textbf{Table 3} (continue)

\noindent The smallest known sizes
$\overline{t}_{2}=\overline{t} _{2}(2,q)<4.75\sqrt{q}$ of
complete arcs in planes $PG(2,q),$

\noindent $2609\leq q\leq 4561$, $A_{q}=\left\lfloor
a_{q}\sqrt{q}-\overline{ t}_{2}(2,q)\right\rfloor $, $B_{q}\geq
\overline{t}_{2}(2,q)/\sqrt{q}$\smallskip

\noindent$\renewcommand{\arraystretch}{1.0}
\begin{array}{@{}r@{\,\,\,\,}c@{\,\,\,\,}c@{\,\,\,}c|@{\,\,}c
@{\,\,\,\,}c@{\,\,\,\,}c@{\,\,\,}c|@{\,\,}c@{\,\,\,\,}
c@{\,\,\,\,}c@{\,\,\,}c|@{\,\,}c@{\,\,\,\,}c@{\,\,\,\,}c@{\,\,\,\,}c@{}}
\hline
q & \overline{t}_{2} & A_{q} & B_{q} & q & \overline{t}_{2} & A_{q} & B_{q}
& q & \overline{t}_{2} & A_{q} & B_{q} & q & \overline{t}_{2} & A_{q} &
B_{q}^{\phantom{H^{L}}} \\ \hline
2903 & 245 & 24 & 4.55 & 3433 & 270 & 22 & 4.61 & 3911 & 292 & 20 & 4.67 &
4423 & 315 & 17 & 4.74 \\
2909 & 246 & 23 & 4.57 & 3449 & 272 & 21 & 4.64 & 3917 & 293 & 19 & 4.69 &
4441 & 315 & 18 & 4.73 \\
2917 & 246 & 24 & 4.56 & 3457 & 271 & 22 & 4.61 & 3919 & 292 & 21 & 4.67 &
4447 & 314 & 19 & 4.71 \\
2927 & 245 & 25 & 4.53 & 3461 & 272 & 22 & 4.63 & 3923 & 293 & 20 & 4.68 &
4451 & 316 & 17 & 4.74 \\
2939 & 246 & 25 & 4.54 & 3463 & 272 & 22 & 4.63 & 3929 & 293 & 20 & 4.68 &
4457 & 315 & 18 & 4.72 \\
2953 & 248 & 23 & 4.57 & 3467 & 272 & 22 & 4.62 & 3931 & 293 & 20 & 4.68 &
4463 & 315 & 19 & 4.72 \\
2957 & 248 & 23 & 4.57 & 3469 & 272 & 22 & 4.62 & 3943 & 294 & 19 & 4.69 &
4481 & 315 & 19 & 4.71 \\
2963 & 248 & 24 & 4.56 & 3481 & 272 & 23 & 4.62 & 3947 & 294 & 20 & 4.68 &
4483 & 317 & 17 & 4.74 \\
2969 & 249 & 23 & 4.57 & 3491 & 273 & 22 & 4.63 & 3967 & 294 & 20 & 4.67 &
4489 & 316 & 19 & 4.72 \\
2971 & 249 & 23 & 4.57 & 3499 & 273 & 22 & 4.62 & 3989 & 296 & 19 & 4.69 &
4493 & 317 & 18 & 4.73 \\
2999 & 250 & 23 & 4.57 & 3511 & 274 & 22 & 4.63 & 4001 & 296 & 20 & 4.68 &
4507 & 318 & 17 & 4.74 \\
3001 & 250 & 23 & 4.57 & 3517 & 274 & 22 & 4.63 & 4003 & 296 & 20 & 4.68 &
4513 & 318 & 17 & 4.74 \\
3011 & 251 & 23 & 4.58 & 3527 & 275 & 21 & 4.64 & 4007 & 297 & 19 & 4.70 &
4517 & 317 & 19 & 4.72 \\
3019 & 251 & 23 & 4.57 & 3529 & 275 & 22 & 4.63 & 4013 & 296 & 20 & 4.68 &
4519 & 319 & 17 & 4.75 \\
3023 & 251 & 23 & 4.57 & 3533 & 275 & 22 & 4.63 & 4019 & 296 & 20 & 4.67 &
4523 & 318 & 18 & 4.73 \\
3037 & 252 & 23 & 4.58 & 3539 & 275 & 22 & 4.63 & 4021 & 296 & 21 & 4.67 &
4547 & 319 & 18 & 4.74 \\
3041 & 252 & 23 & 4.57 & 3541 & 276 & 21 & 4.64 & 4027 & 296 & 21 & 4.67 &
4549 & 319 & 18 & 4.73 \\
3049 & 252 & 24 & 4.57 & 3547 & 276 & 21 & 4.64 & 4049 & 298 & 20 & 4.69 &
4561 & 319 & 18 & 4.73 \\
3061 & 253 & 23 & 4.58 & 3557 & 277 & 21 & 4.65 & 4051 & 299 & 19 & 4.70 &
&  &  &  \\
3067 & 253 & 23 & 4.57 & 3559 & 276 & 22 & 4.63 & 4057 & 298 & 20 & 4.68 &
&  &  &  \\ \hline
\end{array}
$
\end{figure}

From Tables 3 and 4, we obtain Theorem \ref {Th3_5sqroot(q)}.

\begin{thm}
\label{Th3_5sqroot(q)} In $PG(2,q),$ the following holds.
\begin{eqnarray}
t_{2}(2,q)&<&4.75\sqrt{q}~ \mbox{ for }q\leq 4561.
\label{eq3_<4.62sqroot(q)}\\
t_{2}(2,q) &<&4.6\sqrt{q}~~\mbox{ for }q\leq
3209,\,q=3221,3229,3251,3253,3271,3299,3301,3323,\notag\\
&&\phantom{4.6\sqrt{q}~~\mbox{ for }q\leq
3209,\,q=}~3343,3347;\notag\\
t_{2}(2,q) &<&4.7\sqrt{q}~~\mbox{ for }q\leq
4093,\,q=4099,4111,4129,4133,4139,4153,4157,4159,\notag\\
&&\phantom{4.7\sqrt{q}~~\mbox{ for }q\leq
4057,\,q=}~4177,4217,4219,4229,4241,4243,4253,4271,\notag\\
&&\phantom{4.7\sqrt{q}~~\mbox{ for }q\leq
4057,\,q=}~4273,4297,4363.\notag
\end{eqnarray}
  Also,
\begin{eqnarray*}
t_{2}(2,q) &<&5\sqrt{q}-22\mbox{ for }q\leq
3391,\,q= 3413,3433, 3457,3461,3463,3467,3469,3481,\notag \\
&&\phantom{5\sqrt{q}-22\,~ \mbox{ for } q\leq 1369, q=}
3491,3499,3511,3517,3529,3533,3539,3559,\notag
\end{eqnarray*}
\begin{eqnarray*}
&&\phantom{5\sqrt{q}-22\,~ \mbox{ for } q\leq 1369, q=}3583,
3607,3617,3739; \notag\\
t_{2}(2,q) &<&5\sqrt{q}-21\mbox{ for }q\leq 3637,\,
q=3659,3673,3677,3697,3701,3721,3727,3733,\notag\\
&&\phantom{5\sqrt{q}-22\,~ \mbox{ for } q\leq 1369, q=}
3739,3779,3797,3803,3821,3853,3919,4021,\notag\\
&&\phantom{5\sqrt{q}-22\,~ \mbox{ for } q\leq 1369, q=}4027,4153;  \notag\\
t_{2}(2,q) &<&5\sqrt{q}-20\mbox{ for }q\leq 3911,\,q=3919,3923,3929,
3931,3947,3967,4001,4003,\notag\\
&&\phantom{5\sqrt{q}-20\mbox{ for }q\leq 3911,\,q=\,\,}
4013,4019,4021,4027,4049,4057,4073,4099,\notag\\
&&\phantom{5\sqrt{q}-20\mbox{ for }q\leq 3911,\,q=\,\,}
4153,4159,4177,4229,4253,4273,4363;\notag\\
t_{2}(2,q) &<&5\sqrt{q}-19\mbox{ for }q\leq 4297,\,q=
4357,4363,4447,4463,4481,4489,4517.\notag
\end{eqnarray*}
\end{thm}

\begin{figure}[th]
\noindent \textbf{Table 4}

\noindent The smallest known sizes
$\overline{t}_{2}=\overline{t}_{2}(2,q)$ of complete arcs in
planes $PG(2,q)$ with $q\in T_{2}\cup T_{3}$

\noindent $A_{q}=\left\lfloor
5\sqrt{q}-\overline{t}_{2}(2,q)\right\rfloor ,$
$B_{q}>\overline{t}_{2}(2,q)/\sqrt{q}$\smallskip

\noindent$\renewcommand{\arraystretch}{1.0}%
\begin{array}{@{}c@{~\,}c@{~\,}c@{~\,}c@{~\,}c@{\,}|cc@{~\,}c@{~\,}cc@{\,}|c@{~\,}c@{~}c@{~\,}cl@{}}
\hline
q & \overline{t}_{2} & A_{q} & B_{q} & \overline{D}_{q}(\frac{3}{4}) & q &
\overline{t}_{2} & A_{q} & B_{q} & \overline{D}_{q}(\frac{3}{4}) & q &
\overline{t}_{2} & A_{q} & B_{q} & \overline{D}_{q}(\frac{3}{4})^{\phantom{H^{H}}} \\ \hline
4597 & 321 & 18 & 4.74 & 0.9567 & 4813 & 330 & 16 & 4.76 & 0.9573 & 6011 &
377 & 10 & 4.87 & 0.9598 \\
4703 & 325 & 17 & 4.74 & 0.9557 & 4831 & 329 & 18 & 4.74 & 0.9523 & 8192 &
450 & 2 & 4.98 & 0.9560 \\
4723 & 326 & 17 & 4.75 & 0.9563 & 5003 & 338 & 15 & 4.78 & 0.9584 & 2^{14} &
665 &  & 5.20 & 0.9449 \\
4733 & 326 & 17 & 4.74 & 0.9551 & 5347 & 352 & 13 & 4.82 & 0.9599 & 2^{15} &
993 &  & 5.49 & 0.9474 \\
4789 & 329 & 17 & 4.76 & 0.9572 & 5641 & 363 & 12 & 4.84 & 0.9592 & 2^{18} &
3066 &  & 5.99 & 0.9020 \\
4799 & 330 & 16 & 4.77 & 0.9590 & 5843 & 371 & 11 & 4.86 & 0.9604 &  &  &  &
&  \\ \hline
\end{array}
$
\end{figure}

\section{Observations of $\overline{t}_{2}(2,q)$ values}

\label{sec4_observations}

We look for upper estimates of the collection of
$\overline{t}_{2}(2,q)$ values from Tables 1-4 in the
form~(\ref{eq1_KimVu_c=300}), see \cite{KV} and
\cite[Tab.\thinspace 2.6]{HirsStor-2001}. For definiteness, we
use the natural logarithms. Let $c$ be a constant independent
of $q$. We introduce $D_{q}(c)$ and $\overline{D}_{q}(c)$ as
follows:
\begin{eqnarray}
t_{2}(2,q) &=&D_{q}(c)\sqrt{q}\ln ^{c}q,~  \notag \\
\overline{t}_{2}(2,q) &=&\overline{D}_{q}(c)\sqrt{q}\ln ^{c}q.
\label{eq4_Dq(c)}
\end{eqnarray}

Let $\overline{D}_{\text{aver}}(c,q_{0})$ be the average value
of $\overline{ D}_{q}(c)$ calculated in the region $q_{0}\leq
q\leq 4561$ and $q\in T_{2}$ under condition $ q\notin N.$

From Tables 1-4, we obtain Observation~1.\smallskip

\noindent \textbf{Observation 1.} \emph{Let }$173\leq q\leq
4561$ or $q\in T_{2}$, \emph{\ under condition }$q\notin
N$\emph{. Then}

\textbf{(i)} \emph{When $q$ grows,}
$\overline{D}_{q}(0.8)$\emph{\ has a tendency to decreasing.}

(\textbf{ii)} \emph{When $q$ grows,}
$\overline{D}_{q}(0.5)$\emph{\ has a tendency to increasing.}

\textbf{(iii)} \emph{When $q$ grows}, \emph{the values of }$
\overline{D}_{q}(0.75)$\emph{\ oscillate about the average
value }\\$\overline{ D}_{\text{aver}}(0.75,173)=0.95821$ (see
Fig.~1).\emph{ Also,}
\begin{equation}
\begin{array}{cl}
0.946<\overline{D}_{q}(0.75)<0.9647 & \mbox{if }173\leq q\leq 997,\smallskip \\
0.953<\overline{D}_{q}(0.75)<0.9631 & \mbox{if }1009\le q\leq 1999,\smallskip \\
0.951<\overline{D}_{q}(0.75)<0.9618 & \mbox{if }2003\le q\leq 2999,\smallskip \\
0.952<\overline{D}_{q}(0.75)<0.9610 & \mbox{if }3001\le q\leq 3989,\smallskip \\
0.952<\overline{D}_{q}(0.75)<0.9605 & \mbox{if }4001\le q.
\end{array}
\label{eq4_Dq74}
\end{equation}
\begin{figure}[t]
\epsfig{file=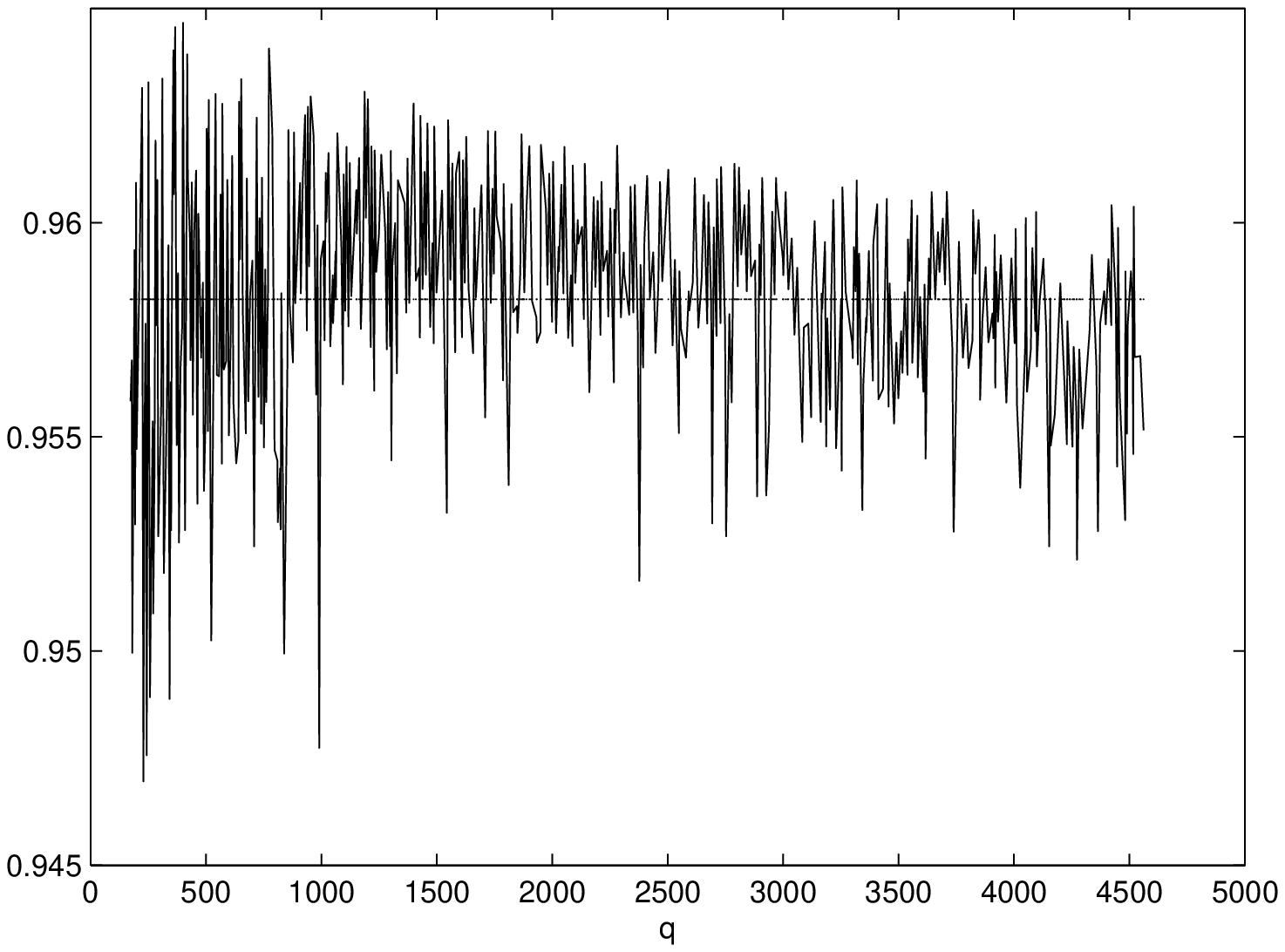,width=\textwidth}
\caption{The values of $\overline{D}_{q}(0.75)$ for $173\le q\le 4561$, $ q\notin N$.
$\overline{ D}_{\text{aver}}(0.75,173)=0.95821$}
\end{figure}
 Data
for relatively big $q$, collected in Table 4, in large confirm
Observation 1.

By Observation 1 it seems that the values of $D_{q}(0.75)$ and
$ \overline{D}_{q}(0.75)$ are sufficiently convenient for
estimates of $ t_{2}(2,q)$ and $\overline{t}_{2}(2,q).$

From Tables 1-4, we obtain Theorem \ref{th4_ln0.75}.

\begin{thm}
\label{th4_ln0.75} In $PG(2,q),$
\begin{equation}
t_{2}(2,q)<0.9987\sqrt{q}\ln ^{0.75}q\quad \mbox{ for }23\leq q\leq 4561,
\text{ }q\in T_{2}\cup T_{3}.  \label{eq4_9987}
\end{equation}
\end{thm}

In Theorem \ref{th1_ln0.75} we slightly rounded the estimate
(\ref{eq4_9987} ).

The graphs of values of $\sqrt{q}\ln ^{0.8}q$, $\sqrt{q}\ln
^{0.75}q$, $ \overline{t}_{2}(2,q)$, and $\sqrt{q}\ln ^{0.5}q$
are shown \smallskip on Fig.~2 where $\sqrt{q}\ln ^{0.8}q$ is
the top curve and $\sqrt{q}\ln ^{0.5}q$
\smallskip is the bottom one.
\begin{figure}[t]
\epsfig{file=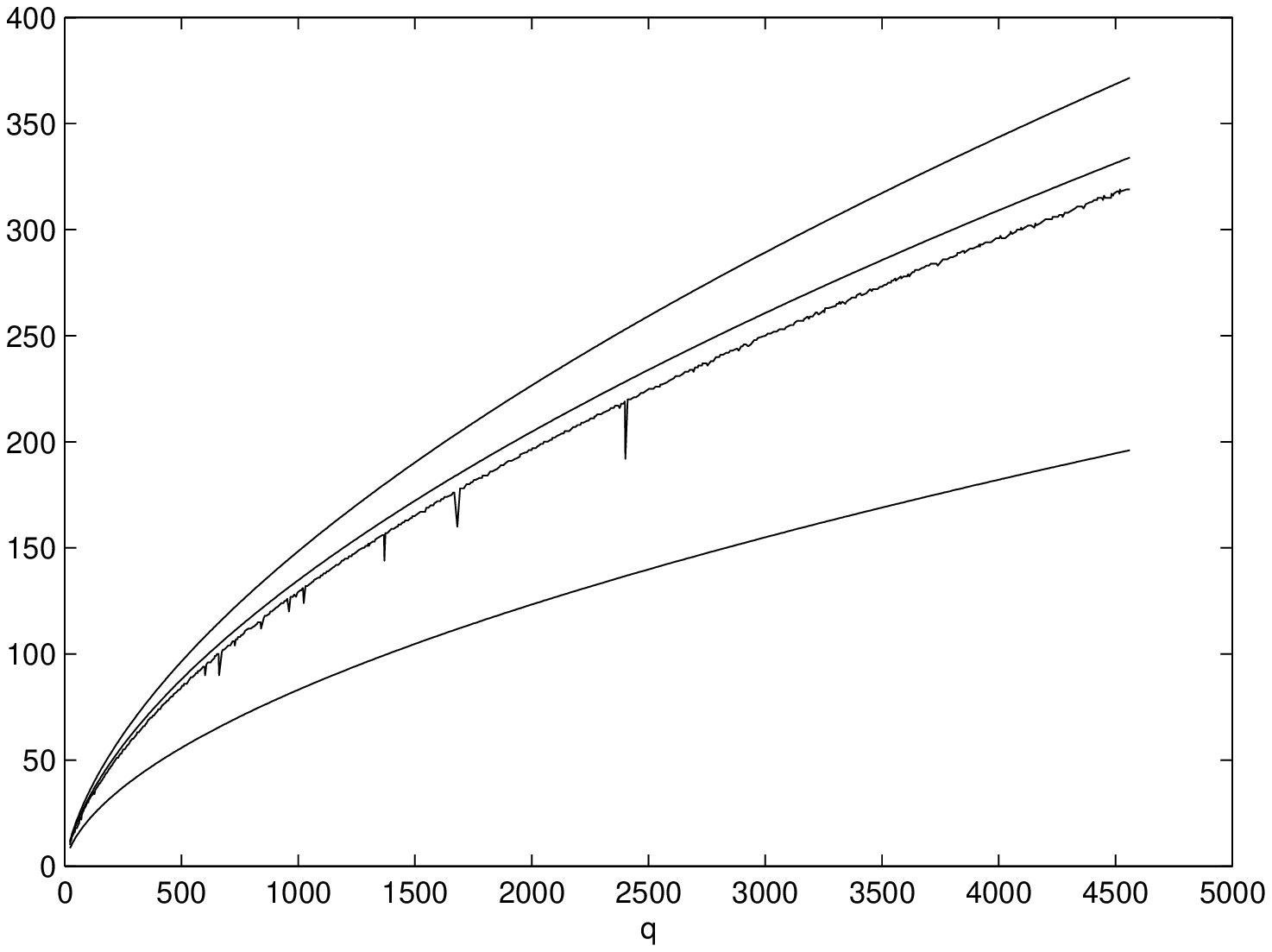,width=\textwidth}
\caption{The values of $\sqrt{q}\ln^{0.8}q$ (the top curve),
$\sqrt{q}\ln^{0.75}q$ (the 2-nd curve), $\overline{t}_{2}(2,q)$ (the 3-rd curve), and
$\sqrt{q}\ln^{0.5}q$ (the bottom curve) for $23\le q\le 4561$}
\end{figure}

One can see on Fig.~2 that always $\overline{t}_{2}(2,q)<\sqrt{q}\ln^{0.75}q$
\smallskip and, moreover, when $q$ grows, the graphs $\sqrt{q}\ln^{0.75}q$
and $\overline{t}_{2}(2,q)$ diverge so that positive difference
$\sqrt{q} \ln^{0.75}q - \overline{t}_{2}(2,q)$ increases.

We denote
\begin{equation}
\widehat{t}_{2}(2,q)=\overline{D}_{\text{aver}}(0.75,173)\sqrt{q}\ln
^{0.75}q,\quad \overline{\Delta }_{q}=\overline{t}_{2}(2,q)-\widehat{t}
_{2}(2,q),\quad \overline{P}_{q}=\frac{100\overline{\Delta }_{q}}{\overline{
t }_{2}(2,q)}\%.  \label{eq4_Delta-q_Pq}
\end{equation}
One can treat $\widehat{t}_{2}(2,q)$ as a \emph{predicted}
value of $ t_{2}(2,q)$. Then $\overline{\Delta }_{q}$ is the
difference between the smallest known size
$\overline{t}_{2}(2,q)$ of complete arcs and the predicted
value. Finally, $\overline{P}_{q}$ is this difference in
percentage terms of the smallest known size.\smallskip

\noindent \textbf{Observation 2.} \emph{Let }$173\leq q\leq
4561$\emph{\ or }$ q\in T_{2}$, $ q\notin N$\emph{.
Then\newline }
\begin{equation}
-2.05<\overline{\Delta }_{q}<0.83.  \label{eq4_Delta_q}
\end{equation}
\begin{equation}
\begin{array}{cl}
-1.19\%<\overline{P}_{q}<0.67\% & \mbox{if }131\leq q\leq 997,\smallskip  \\
-0.53\%<\overline{P}_{q}<0.51\% & \mbox{if }1009\le q\leq 1999,\smallskip  \\
-0.70\%<\overline{P}_{q}<0.38\% & \mbox{if }2003\le q\leq 2999,\smallskip  \\
-0.57\%<\overline{P}_{q}<0.29\% & \mbox{if }3001\le q\leq 3989,\smallskip  \\
-0.64\%<\overline{P}_{q}<0.23\% & \mbox{if }4001\le q.
\end{array}
\label{eq4_percent}
\end{equation}

 By
(\ref{eq4_Delta_q}) and (\ref{eq4_percent}), see also Fig. 3
and 4, the upper bounds of $\overline{\Delta }_{q}$ and
$\overline{P}_{q}$ are relatively small. Moreover, the upper
bound of $\overline{P}_{q}$ decreases when $q$ grows. Therefore
the values of $\overline{\Delta }_{q}$ and $ \overline{P}_{q}$
are useful for computer search of small arcs.
\begin{figure}[tb]
\epsfig{file=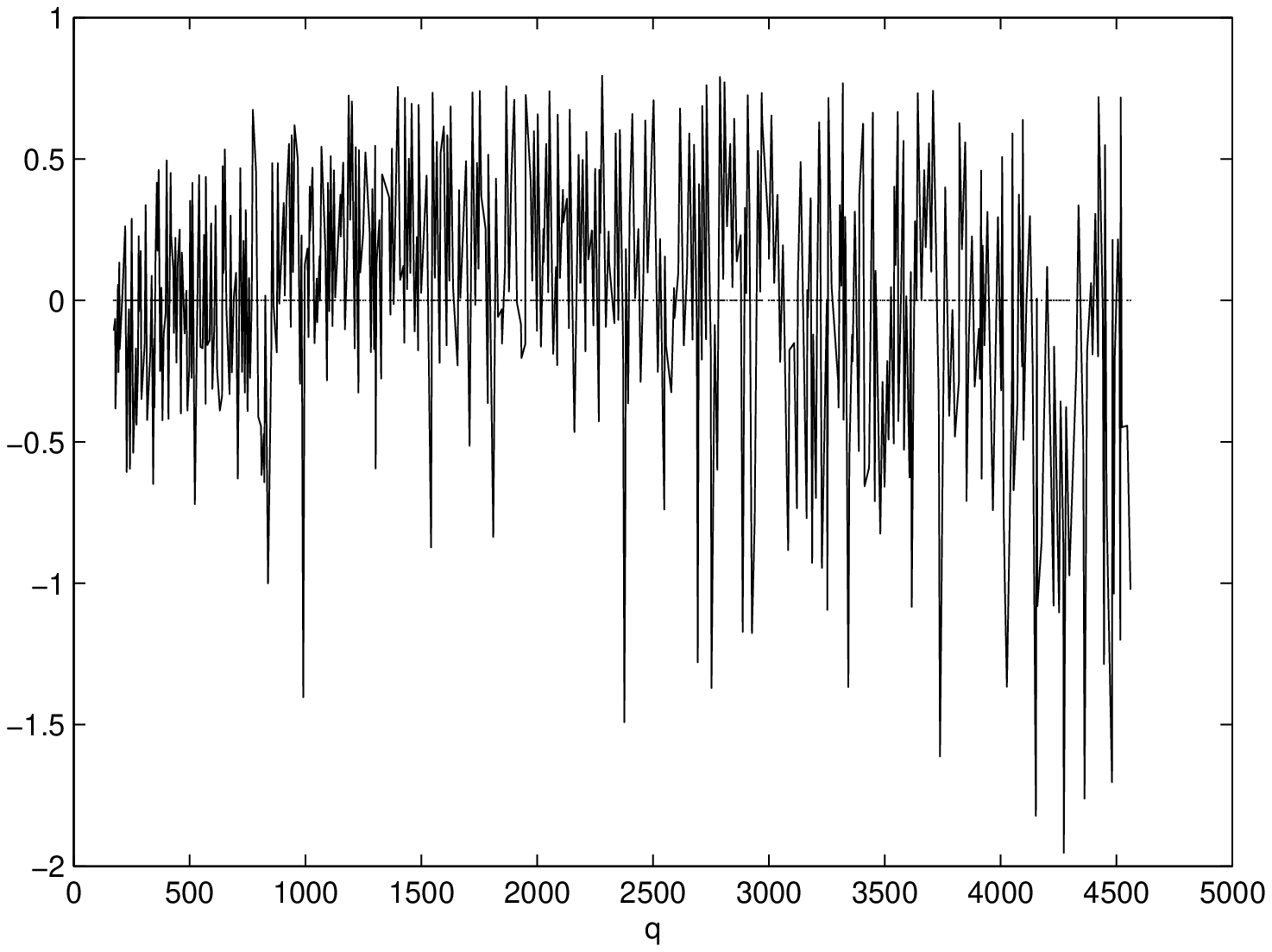,width=\textwidth}
\caption{The values of $\overline{\Delta }_{q}=\overline{t}_{2}(2,q)-\widehat{t}
_{2}(2,q)$ for $173\le q\le 4561$,  $ q\notin N$}
\end{figure}
\begin{figure}[tb]
\epsfig{file=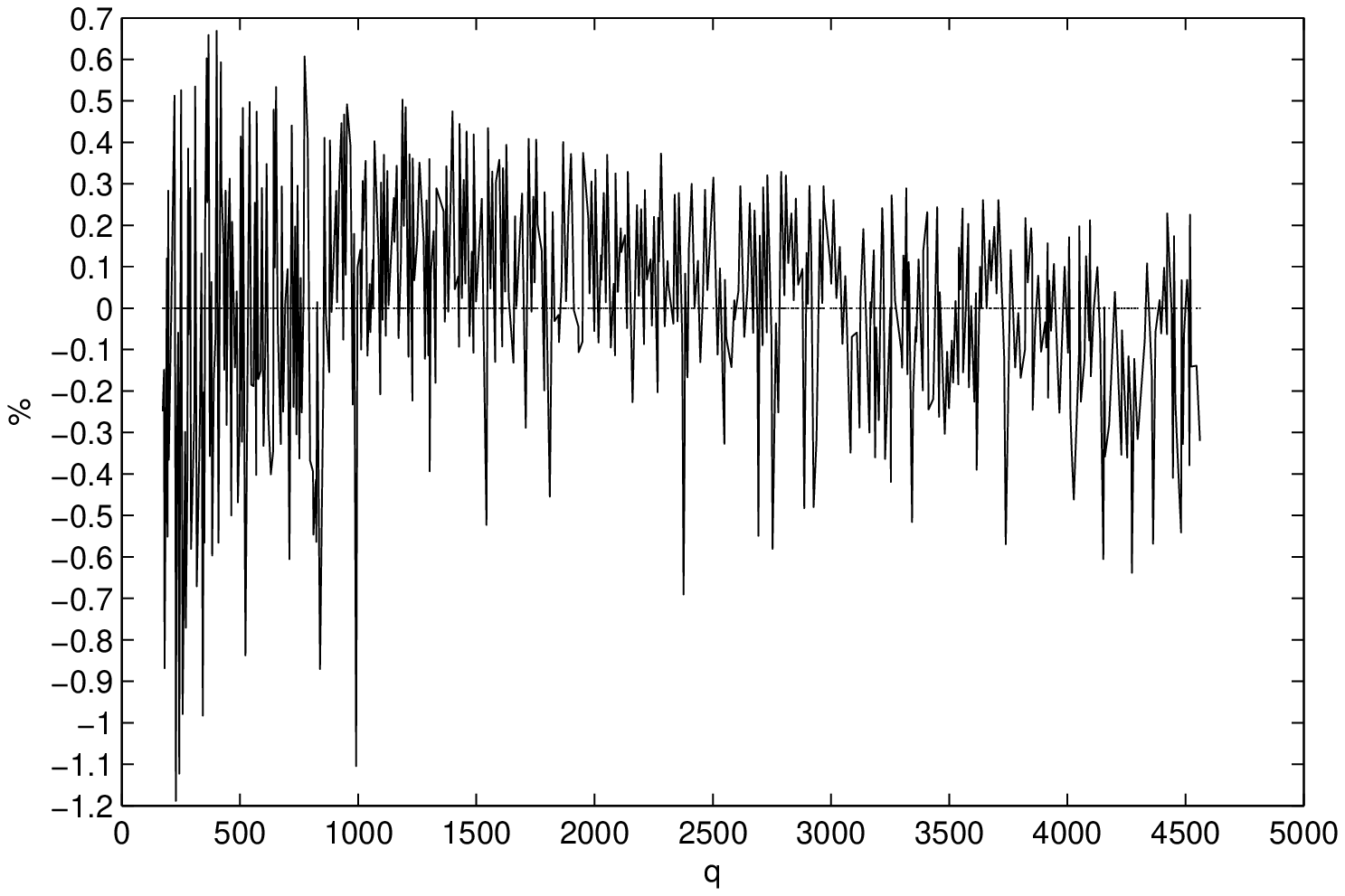,width=\textwidth}
\caption{The values of $\overline{P}_{q}=100\overline{\Delta }_{q}/\overline{
t }_{2}(2,q)\%$ for $173\le q\le 4561$,  $ q\notin N$}
\end{figure}

The relations (\ref{eq4_Dq74})--(\ref{eq4_percent}), Theorems
\ref{th1_ln0.75} and \ref{th4_ln0.75}, and Figures 1--4 are the
foundation for Conjecture \ref {conj1_ln0.75}.

\begin{rmk}
\label{rem4} By above, $\sqrt{q}\ln^{0.75}q$ seems a reasonable
upper bound on the current collection of
$\overline{t}_{2}(2,q)$ values. It gives some reference points
for computer search and foundations for Conjecture \ref
{conj1_ln0.75} on the upper bound for $t_{2}(2,q)$. In
principle, the constant $c=0.75$ can be sightly reduced to move
the curve $\sqrt{q}\ln^{c}q$ near to the curve of
$\overline{t}_{2}(2,q)$, see Fig.~2. For example, from Tables
1-4, Theorem \ref{th4_ln0.741} holds.
\end{rmk}

\begin{thm}
\label{th4_ln0.741} In $PG(2,q),$
\begin{equation}
t_{2}(2,q)<\sqrt{q}\ln ^{0.735}q\quad \mbox{ for }32\leq q\leq 4561,
\text{ }q\in T_{2}\cup T_{3}.  \label{eq4_999_741}
\end{equation}
\end{thm}

\section{Constructions of families of complete arcs in $PG(2,q)$}
\label{sec5_constr} In the homogenous coordinates of a point
$(x_{0},x_{1},x_{2})$ we put $x_{0}\in \{0,1\},$
$x_{1},x_{2}\in F_{q}.$ Let $F_{q}^{\ast }=F_{q}\backslash
\{0\}.$ Let $\xi $ be a primitive element of $F_{q}.$ Remind
that \emph{indexes of powers of }$\xi $\emph{\ are calculated
modulo }$q-1.$

Throughout this section we use the conic $\mathcal{C}$ of
equation $x_{1}^{2}=x_{0}x_{2}.$ We denote points of
$\mathcal{C}$ as follows:
\begin{equation*}
A_{i}=(1,i,i^{2}),\text{ }i\in F_{q};\text{\quad }\overline{A}_{d}=(1,\xi
^{d},\xi ^{2d}),\text{ }d\in \{0,1,\ldots ,q-2\};\text{\quad }A_{\infty
}=(0,0,1).
\end{equation*}

\subsection{Arcs with two points on a tangent to a conic}

Through this subsection, $q\geq 19$ is an \emph{odd}
\emph{prime}. Let $H$ be an integer in the
region\begin{equation} \left\lfloor \frac{q-1}{3}\right\rfloor
\leq H\leq \frac{q-1}{2}. \label{eq5.1_H}
\end{equation}We denote by $\mathcal{V}_{H}$ the following $(H+1)$-subset of the conic $\mathcal{C}:$
\begin{equation}
\mathcal{V}_{H}=\{A_{i}:i=0,1,2,\ldots ,H\}\subset \mathcal{C}.
\label{eq5.1_VH}
\end{equation}We denote the points of $PG(2,q):$\begin{equation}
P=(0,1,0),\text{ }T_{H}=(0,1,b_{H}),\text{ }b_{H}=\left\{
\begin{array}{ll}
2H+1 & \text{if }H=\left\lfloor \frac{1}{3}(q-1)\right\rfloor \smallskip \\
2H & \text{if }\left\lfloor \frac{1}{3}(q-1)\right\rfloor <H\leq \frac{1}{2}(q-1)\end{array}\right. .  \label{eq5.1_PT_M}
\end{equation}Let $\ell _{0}$ be the line of equation $x_{0}=0.$ It is the
\emph{tangent} to $\mathcal{C}$ at $A_{\infty }.$ It holds that
$\{P$,$T_{H}\}\subset \ell _{0}$.

\noindent \textbf{Construction A. }Let $q$ be an \emph{odd}
\emph{prime}. Let $H,$ $\mathcal{V}_{H},$ $P$ and $T_{H}$ be
given by (\ref{eq5.1_H})--(\ref{eq5.1_PT_M}). We construct a
point $(H+3)$-set $\mathcal{K}_{H}$ in the plane $PG(2,q)$ as
follows:\begin{equation*} \mathcal{K}_{H}=\mathcal{V}_{H}\cup
\{P,T_{H}\}.
\end{equation*}

The following lemma can be proved by elementary calculations.

\begin{lem}
\label{lem5.1_01b_i+j}\emph{(i)}\textbf{\ }Let $i\neq j.$ A
point $(0,1,b)$ is collinear with points $A_{i},$ $A_{j}$ if
and only if
\begin{equation}
b=i+j.  \label{eq5.1_b=i+j}
\end{equation}

\emph{(ii) }Let $i\neq j,$ $a,b\in F_{q},$ $b\neq a^{2}.$ Then
a point $(1,a,b)$ is collinear with $A_{i},$ $A_{j}$ if and
only if $b=a(i+j)-ij.$

\emph{(iii) }Let $a\in F_{q},$ $a\neq i.$ Then a point
$(1,a,i^{2})$ is collinear with $P,~A_{i}.$
\end{lem}

\begin{thm}
The $(H+3)$-set $\mathcal{K}_{H}$ of Construction A is an arc
in $PG(2,q)$.
\end{thm}

\begin{pf}
By (\ref{eq5.1_H}),(\ref{eq5.1_VH}), the sum $i+j$ in
(\ref{eq5.1_b=i+j}) is running on $\{1,2,\ldots ,2H-1\}$ where
$2H-1\leq q-2$ if $\left\lfloor \frac{1}{3}(q-1)\right\rfloor
<H$ and $2H-1<\frac{2}{3}(q-1)$ if $H=\left\lfloor
\frac{1}{3}(q-1)\right\rfloor .$ So, $\{0,b_{H}\}\cap
\{1,2,\ldots ,2H-1\}=\emptyset ,$ see (\ref{eq5.1_PT_M}).
Therefore $P$ and $T_{H}$ do not lie on bisecants
of$~\mathcal{V}_{H}.$ In other side, any point of
$\mathcal{V}_{H}$ does not lie on the line $PT_{H}$ as $PT_{H}$
is a tangent to $\mathcal{C}$ in $A_{\infty }.$ $\hfill \qed$
\end{pf}

\begin{thm}
\label{th5.1_l0||C-VM_onBisecKM}Let $H$ be given by
\emph{(\ref{eq5.1_H})}. Then all points of $\ell _{0}\cup
\mathcal{C}\backslash \mathcal{V}_{H}$ lie on bisecants of
$\mathcal{K}_{H}.$
\end{thm}

\begin{pf}
All points of $\ell _{0}$ are covered as two points $P$ and
$T_{H}$ of this line belong to $\mathcal{K}_{H}.$

Let $\mathcal{R}$ and $\mathcal{S}$ be sets of integers modulo
$q,$ i.e. $\mathcal{R}\cup \mathcal{S}\subset F_{q}$.

Let $\mathcal{R}=\{-H,-(H-1),\ldots ,-1\}=\{q-H,q-(H-1),\ldots
,q-1\}.$ By Lemma~\ref{lem5.1_01b_i+j}(i), points
$A_{j},A_{-j},P$ are collinear. Therefore, a point $A_{j}$ of
$\mathcal{C}\backslash \mathcal{V}_{H}$ with $j\in \mathcal{R}$
lies on the bisecant of $\mathcal{K}_{H}$ through $P$ and
$A_{-j}$ where $-j\in \{H,H-1,\ldots ,1\},$ $A_{-j}\in
\mathcal{V}_{H}.$

Let $\mathcal{S}=\{b_{H}-H,b_{H}-(H-1),\ldots
,b_{H}-1,b_{H}\}.$ By Lemma \ref{lem5.1_01b_i+j}(i), a point
$A_{j}\ $of $\mathcal{C}\backslash \mathcal{V}_{H}$ with $j\in
\mathcal{S}$ lies on the bisecant $T_{H}A_{b_{H}-j}$ where
$b_{H}-j\in \{H,H-1,\ldots ,1,0\},$ $A_{b_{H}-j}\in
\mathcal{V}_{H}.$

Let $b_{H}=2H+1.$ Then $\mathcal{S}=\{H+1,H+2,\ldots ,2H+1\}.$
Also, by (\ref{eq5.1_PT_M}), $H=\left\lfloor
\frac{1}{3}(q-1)\right\rfloor $ whence $H= \frac{1}{3}(q-v)$
where $v\in\{1,2\}$ and $v\equiv q ~(\bmod~3)$. Hence,
$3H=q-v$, $2H+1=q-H-v+1\in\{q-H,q-H-1\}.$

Let $b_{H}=2H.$ Then $\mathcal{S}=\{H,H+1,\ldots ,2H\}.$ Also,
by (\ref{eq5.1_PT_M}), $H>\left\lfloor
\frac{1}{3}(q-1)\right\rfloor $ whence $H> \frac{1}{3}(q-v)$
where $v\in\{1,2\}$ is as above. Therefore $3H> q-v$, $2H>
q-H-v$, $2H\ge q-H-v+1\in\{q-H,q-H-1\}.$

We proved that $\{H+1,H+2,\ldots ,q-1\}\subseteq
\mathcal{S}\cup \mathcal{R} . $ Also we showed that the points
$A_{j}\ $of $\mathcal{C}\backslash \mathcal{V}_{H}$ with $j\in
\mathcal{S}\cup \mathcal{R}$ are covered by bisecants of
$\mathcal{K}_{H}$ through $P$ (if $j\in \mathcal{R}$) or
through $T_{H}$ (if $j\in \mathcal{S)}.$ In the other side,
$\mathcal{C} \backslash
\mathcal{V}_{H}=\{A_{j}:j=H+1,H+2,\ldots ,q-1\}\cup \{A_{\infty
}\}$ where $A_{\infty }\in \ell _{0}$. So, all points of
$\mathcal{C} \backslash \mathcal{V}_{H}$ are covered$.$ $\hfill
\qed$
\end{pf}

\begin{definition}
\label{def5.1_overlineH}Let $q$ be an \emph{odd} \emph{prime}.
Let $ \overline{H}$ be an integer and let
\begin{equation*}
\mathcal{P}_{\overline{H}}=\{P\}\cup \{A_{i}:i=0,1,2,\ldots ,\overline{H}\}.
\end{equation*}
We call \emph{critical value} \emph{of }$\overline{H}$ and
denote by $ \overline{H}_{q}$ the \emph{smallest} value of
$\overline{H}$ such that all points of the form $(1,a,b),$
$a,b\in F_{q},$ $b\neq a^{2},$ lie on bisecants of
$\mathcal{P}_{\overline{H}}.$
\end{definition}

\begin{thm}
Let $q\geq 19$ be an \emph{odd} \emph{prime}. Let
$\overline{H}_{q}\leq \frac{1}{2}(q-1)$ and let
\begin{equation*}
\max \{\overline{H}_{q},\left\lfloor \frac{q-1}{3}\right\rfloor \}\leq H\leq
\frac{q-1}{2}.
\end{equation*}
Then the arc $\mathcal{K}_{H}$ of Construction A is complete.
\end{thm}

\begin{pf}
We use Theorem \ref{th5.1_l0||C-VM_onBisecKM} and Definition
\ref {def5.1_overlineH}.$\hfill \qed$
\end{pf}

In this subsection, we put $q\geq 19$ as we checked by computer
that $\frac{1 }{2}(q-1)<\overline{H}_{q}$ if $q\leq 17.$

\begin{corollary}
\label{cor5.1_familyFq}Let $q\geq 19$ be an \emph{odd}
\emph{prime}. Let $ \overline{H}_{q}\leq \frac{1}{2}(q-1).$
Then Construction A forms a family of complete $k$-arcs in
$PG(2,q)$ containing arcs of all sizes $k$ in the region
\begin{equation*}
\max \{\overline{H}_{q},\left\lfloor \frac{q-1}{3}\right\rfloor \}+3\leq
k\leq \frac{q+5}{2}.
\end{equation*}
If $\overline{H}_{q}\leq \left\lfloor
\frac{1}{3}(q-1)\right\rfloor $ then cardinality of this family
is equal to $\left\lceil \frac{1}{6} (q+5)\right\rceil $ and
size of the smallest complete arc of the family is $
\left\lfloor \frac{1}{3}(q+8)\right\rfloor .$
\end{corollary}

By computer search using Lemma \ref{lem5.1_01b_i+j}(ii),(iii)
we obtained the following theorem.

\begin{thm}
\label{th5.1_byTable1}Let $q\geq 19$ be an \emph{odd}
\emph{prime}.
\smallskip Let $\overline{H}_{q}$ be\ given by Definition \emph{\ref
{def5.1_overlineH}}. We introduce $D_{q}$ and $\Delta _{q}$ as
follows: $ \overline{H}_{q}=D_{q}\sqrt{q}\ln ^{0.9}q,$ $\Delta
_{q}=\left\lfloor \frac{1 }{3}(q-1)\right\rfloor
-\overline{H}_{q}.$ Then the following holds.
\begin{eqnarray}
\emph{(i)\ }&&\left\lfloor \frac{q-1}{3}\right\rfloor <\overline{H}_{q}\leq
\frac{q-1}{2}\text{ if }19\leq q\leq 71\text{ and }q=79,83,89,107.\notag\\
\text{\emph{(ii) }}&&\overline{H}_{q}\leq \left\lfloor
\frac{q-1}{3} \right\rfloor \text{if }109\leq
q\leq 1367,~2003\leq
q\leq 2063,~q=73,97,101,103.\notag\\
\text{\emph{(iii) }}&&\overline{H}_{q}<0.98\sqrt{q}\ln ^{0.9}q~\text{ if }
19\leq q\leq 1367,\,2003\leq
q\leq 2063.  \label{eq5.1_ln^0.9}\\
\text{\emph{(iv) }}&&
\begin{array}{cll}
0\leq \Delta _{q}\leq 99, & 0.74<D_{q}<0.98, & \text{if }109\leq q\leq 599;
\\
89\leq \Delta _{q}\leq 198, & 0.75<D_{q}<0.94, & \text{if }601\leq q\leq
1049; \\
187\leq \Delta _{q}\leq 288, & 0.72<D_{q}<0.95, & \text{if }1051\leq q\leq
1367;\\
417\leq \Delta _{q}\leq 464, &0.78 <D_{q}<0.91, & \text{if }2003\leq q\leq
2063 .
\end{array}  \notag
\end{eqnarray}
\end{thm}

For situations $\left\lfloor \frac{1}{3}(q-1)\right\rfloor
<\overline{H} _{q}, $ the values of $\overline{H}_{q}$ are
given in Table 5.

\begin{figure}[th]
\noindent \textbf{Table 5}

\noindent The values
$\overline{H}_{q},\overline{G}_{q},\overline{J}_{q}$ for cases
$\overline{H}_{q},\overline{G}_{q}>\left\lfloor \frac{1}{3}
(q-1)\right\rfloor
,\overline{J}_{q}>\frac{1}{4}(q-3)$\smallskip

\noindent $\renewcommand{\arraystretch}{1.0}
\begin{array}{@{}r@{~\,}rrr|r@{~\,}rrr|r@{~\,}rrr|r@{~\,}rrr|rc@{}}
\hline
q & \overline{H}_{q} & \overline{G}_{q} & \overline{J}_{q} & q & \overline{H}
_{q} & \overline{G}_{q} & \overline{J}_{q} & q & \overline{H}_{q} &
\overline{G}_{q} & \overline{J}_{q} & q & \overline{H}_{q} & \overline{G}_{q}
& \overline{J}_{q} & q & \overline{J}_{q}^{\phantom{{H}}} \\ \hline
19 & 9 &  &  & 43 & 19 &  & 16 & 71 & 24 &  & 25 & 103 &  &  & 34 & 167 & 42
\\
23 & 11 &  &  & 47 & 18 &  & 18 & 73 &  & 27 &  & 107 & 36 &  & 30 & 179 & 62
\\
27 &  &  & 12 & 49 &  & 18 &  & 79 & 29 &  & 27 & 125 &  & 43 &  & 191 & 49
\\
29 & 13 &  &  & 53 & 20 & 19 &  & 81 &  & 29 &  & 127 &  &  & 39 & 211 & 54
\\
31 & 14 &  & 14 & 59 & 23 &  & 24 & 83 & 29 &  & 29 & 131 &  &  & 38 & 223 &
60 \\
32 &  & 15 &  & 61 & 22 & 22 &  & 89 & 32 &  &  & 139 &  &  & 38 & 343 & 86
\\
37 & 16 & 16 &  & 64 &  & 24 &  & 97 &  & 33 &  & 151 &  &  & 42 &  &  \\
41 & 16 & 19 &  & 67 & 24 &  & 23 & 101 &  & 35 &  & 163 &  &  & 41 &  &  \\
\hline
\end{array}
$
\end{figure}

\begin{thm}\label{th5.1_result}
Let $q$ be an \emph{odd} \emph{prime} with $109\leq q\leq
1367$, $2003\le q\le 2063$, or $ q=73,97,101,103.$ Then
Construction A forms a family of complete $k$-arcs in $PG(2,q)$
containing arcs of all sizes $k$ in the region
\begin{equation*}
\left\lfloor \frac{q+8}{3}\right\rfloor \leq k\leq \frac{q+5}{2}.
\end{equation*}
\end{thm}

\begin{pf}
We use Corollary \ref{cor5.1_familyFq} and Theorem
\ref{th5.1_byTable1}(ii).$ \hfill \qed$
\end{pf}

\subsection{Arcs with two points on a bisecant of a conic}

Throughout this subsection, $q\geq 32$ is a \emph{prime power}.
Let $G$ be an integer in the region
\begin{equation}
\left\lfloor \frac{q-1}{3}\right\rfloor \leq G\leq \left\lceil \frac{q-3}{2}
\right\rceil .  \label{eq5.2_G}
\end{equation}
We denote by $\mathcal{D}_{G}$ the following $(G+1)$-subset of
the conic $ \mathcal{C}:$
\begin{equation}
\mathcal{D}_{G}=\{\overline{A}_{d}:d=0,1,2,\ldots ,G\}\subset \mathcal{C}.
\label{eq5.2_DG}
\end{equation}
Clearly, $A_{0}\notin \mathcal{D}_{G}.$ Let
\begin{equation}
\gamma \in F_{q},\text{ }\gamma =\left\{
\begin{array}{cl}
-1=\xi ^{(q-1)/2} & \text{for }q\text{ odd} \\
1 & \text{for }q\text{ even}
\end{array}
\right. .  \label{eq5.2_gamma_-1_1}
\end{equation}
We denote the points of $PG(2,q):$
\begin{equation}
Z=(1,0,\gamma \xi ^{0}),\text{ }B_{G}=(1,0,\gamma \xi ^{\beta _{G}}),\text{ }
\beta _{G}=2G.  \label{eq5.2_Z_BG}
\end{equation}
Let $\ell _{1}$ be the line of equation $x_{1}=0.$ It is the
\emph{bisecant} $A_{\infty }A_{0}$ of $\mathcal{C}.$ We have
$\{Z,B_{G}\}\subset \ell _{1}$.

Using (\ref{eq5.2_G}),(\ref{eq5.2_gamma_-1_1}), by elementary
calculations we obtained the following lemma.

\begin{lem}
\label{lem5.2_(0,1,b)}\emph{(i) Let }$d\neq t.$ A point
$(0,1,b)$ is collinear with points $\overline{A}_{d},$
$\overline{A}_{t}$ if and only if
\begin{equation}
b=\xi ^{d}+\xi ^{t}.  \label{eq5.2_b=xi^i+xi^j}
\end{equation}

\emph{(ii) }A point $(0,1,b)$ is collinear with points
$\overline{A}_{d}$ and $(1,0,\gamma U)$ if and only if
\begin{equation}
b=\frac{\xi ^{2d}+U}{\xi ^{d}}.  \label{eq5.2_b=xi^2i+U}
\end{equation}
\end{lem}

\begin{corollary}
\label{cor5.2_(0,1,0)}\emph{(i) }For all $q,$ the point
$P=(0,1,0)$ does not lie on any bisecant of $\mathcal{D}_{G}$.

\emph{(ii) }Let $q$ be even. Then the points
$P,Z,\overline{A}_{d}$ are collinear if and only if $d=0.$
Also, the points $P,B_{G},\overline{A}_{d}$ are collinear if
and only if $d=G$.

\emph{(iii) }Let $q\equiv 1$ $(\bmod~4).$ Then the points $P,Z,
\overline{A}_{d}$ are collinear if and only if $d\in
\{\frac{1}{4}(q-1),
\frac{3}{4}(q-1)\}.$ Also, the points
$P,B_{G},\overline{A}_{d}$ are collinear if and only if
$\smallskip d\in \{G+\frac{1}{4}(q-1),$ $G-\frac{1}{
4}(q-1)\}$.

\emph{(iv) }Let $q\equiv 3$ $(\bmod~4).$ Then the points $P,Z,
\overline{A}_{d}$ and $P,B_{G},\overline{A}_{d}$ are not
collinear for any$ ~d.$
\end{corollary}

\begin{pf}
(i) In (\ref{eq5.2_b=xi^i+xi^j})$,$ the case $b=0$ implies $\xi
^{d}+\xi ^{t}=0.$ For even $q,$ it is impossible$.$ For odd
$q,$ we obtain $\xi ^{d}=-\xi ^{t}$ whence, by
(\ref{eq5.2_gamma_-1_1}), $d=t+(q-1)/2.$ By (\ref {eq5.2_G}),
it is impossible.

(ii)-(iv) In (\ref{eq5.2_b=xi^2i+U})$,$ the case $b=0$ implies
$\xi ^{2d}+U=0 $ whence $\xi ^{2d}+1=0$ if $(1,0,\gamma U)=Z$
and $\xi ^{2d}+\xi ^{\beta _{G}}=0$ if $(1,0,\gamma U)=B_{G}.$
Remind that $\beta _{G}=2G$ and $ d\leq q-2.$

(ii) Here $q$ is even but $q-1$ is odd. For $(1,0,\gamma U)=Z,$
we have $\xi ^{2d}=1$ whence $d=0.$ For $(1,0,\gamma U)=B_{G},$
it holds that $\xi ^{2d}=\xi ^{\beta _{G}}$ whence $2d\equiv
2G$ $(\bmod q-1).$ So, $ d=G. $

(iii) Here both $q-1$ and $\frac{1}{2}(q-1)$ are even. If
$(1,0,\gamma U)=Z\ $then $\xi ^{2d}=-1=\xi ^{(q-1)/2}$ whence
$2d\equiv \frac{1}{2}(q-1)$ $( \bmod q-1).$ It is possible if
$d=\frac{1}{4}(q-1)$ or $d=\frac{3}{4} (q-1)$. If $(1,0,\gamma
U)=B_{G}\ $then, by (\ref{eq5.2_gamma_-1_1}), $\xi ^{2d}=-\xi
^{\beta _{G}}=\xi ^{\beta _{G}+(q-1)/2}$ whence $2d\equiv
2G+(q-1)/2$ $(\bmod q-1).$ So, $d=G+\frac{1}{4}(q-1)$ or
$d=G-\frac{1 }{4}(q-1)$.

(iv) Here $q-1$ is even whereas $\frac{1}{2}(q-1)$ is odd. If
$(1,0,\gamma U)=Z$ then $\xi ^{2d}=-1=\xi ^{(q-1)/2}$ whence
$2d\equiv \frac{1}{2}(q-1)$ $ (\bmod q-1).$ It is impossible.
For $(1,0,\gamma U)=B_{G},$ it holds that $\xi ^{2d}=-\xi
^{\beta _{G}}=\xi ^{2G+(q-1)/2}$ that is impossible. $ \hfill
\qed$
\end{pf}

\noindent \textbf{Construction B. }Let $q$ be a \emph{prime
power. }Assume that $q\not\equiv 3$ $(\bmod~4).$\emph{\ }Let
$G,$ $\mathcal{D}_{G},$ $Z,$ and $B_{G}$ be given by
(\ref{eq5.2_G})--(\ref{eq5.2_Z_BG}). We construct a point
$(G+3)$-set $\mathcal{W}_{G}$ in $PG(2,q)$ as follows:
\begin{equation*}
\mathcal{W}_{G}=\mathcal{D}_{G}\cup \{Z,B_{G}\}.
\end{equation*}

From Lemma \ref{lem5.1_01b_i+j} it follows.

\begin{lem}
\label{lem5.2_01b_i+j}Let $d\neq t.$ A point $(1,0,\gamma \xi
^{\beta })$ is collinear with $\overline{A}_{d}$,
$\overline{A}_{t}$ if and only if
\begin{equation}
\beta =d+t.  \label{eq5.2_beta=i+j}
\end{equation}
\end{lem}

\begin{thm}
\label{th5.2_is-arc}The $(G+3)$-set $\mathcal{W}_{G}$ of
Construction B is an arc.
\end{thm}

\begin{pf}
By (\ref{eq5.2_G}),(\ref{eq5.2_DG}), the sum $d+t$ in
(\ref{eq5.2_beta=i+j}) is running on $\{1,2,\ldots ,2G-1\}$
where $2G-1\leq q-3.$ So, $\{0,\beta _{G}\}\cap \{1,2,\ldots
,2G-1\}=\emptyset ,$ see (\ref{eq5.2_Z_BG}). Therefore $Z$ and
$B_{G}$ do not lie on bisecants of$~\mathcal{D}_{G}.$ In other
side, any point of $\mathcal{D}_{G}$ does not lie on the line
$ZB_{G}$ as $ZB_{G}$ is the bisecant of $\mathcal{C}$ through
$A_{\infty }$ and $ A_{0} $ where $\{A_{\infty },A_{0}\}\cap
\mathcal{D}_{G}=\emptyset $. $ \hfill \qed$
\end{pf}

\begin{thm}
\label{th5.2_l0||C-VM_onBisecKM}Let $q$ be a \emph{prime power.
}Assume that $q\not\equiv 3$ $(\bmod~4).$\emph{\ }Let $G$ be
given by~\emph{(\ref {eq5.2_G})}. Then all points of $\{P\}\cup
\ell _{1}\cup \mathcal{C} \smallsetminus \mathcal{D}_{G}$ lie
on bisecants of the arc $\mathcal{W}_{G}$ of Construction$~$B.
\end{thm}

\begin{pf}
By (\ref{eq5.2_G}),(\ref{eq5.2_DG}),
$\{\overline{A}_{0},\overline{A} _{(q-1)/4}\}\subset
\mathcal{D}_{G}.$ So, the point $P$ is covered
by Corollary~\ref{cor5.2_(0,1,0)}(ii),(iii).

All points of $\ell _{1}$ are covered as two points $Z$ and
$B_{G}$ of this line belong to $\mathcal{W}_{G}.\smallskip $

Throughout this proof, $\mathcal{R}$ and $\mathcal{S}$ are sets
of integers modulo $q-1.$ It can be said that $\mathcal{R}$ and
$\mathcal{S}$ are sets of indexes of powers of $\xi .$

Let $\mathcal{R}=\{-G,-(G-1),\ldots
,-1\}=\{q-1-G,q-1-(G-1),\ldots ,q-2\}.$ By\smallskip\
Lemma~\ref{lem5.2_01b_i+j}, points $\overline{A}_{t},\overline{
A}_{-t},Z$ are collinear. Therefore, a point $\overline{A}_{t}$
of $\mathcal{ C}\smallsetminus \mathcal{D}_{G}$ with $t\in
\mathcal{R}$ lies on the
\smallskip bisecant of $\mathcal{W}_{G}$ through $\overline{A}_{-t}$ and $Z$
where $-t\in \{G,G-1,\ldots ,1\},$ $\overline{A}_{-t}\in
\mathcal{D}_{G}.$

Let $\mathcal{S}=\{\beta _{G}-G,\beta _{G}-(G-1),\ldots ,\beta
_{G}-1,\beta _{G}\}.$ By Lemma \ref{lem5.2_01b_i+j},\ a point
$\overline{A}_{t}\ $of $ \mathcal{C}\smallsetminus
\mathcal{D}_{G}$ with $t\in \mathcal{S}$ lies on the bisecant
$B_{G}\overline{A}_{\beta _{G}-t}$ where $\beta _{G}-t\in
\{G,G-1,G-2,\ldots ,1,0\}$, $\overline{A}_{\beta _{G}-t}\in
\mathcal{D}_{G}.$

As $\beta _{G}=2G,$ we have $\mathcal{S}=\{G,G+1,\ldots ,2G\}.$
Also, $G\geq \left\lfloor \frac{1}{3}(q-1)\right\rfloor .$ If
$3|(q-1)$ then $G\geq \frac{ 1}{3}(q-1)$ whence $2G\geq q-G-1.$
If $3\nmid (q-1)$ then $G\geq \frac{1}{3} (q-2)$ whence $2G\geq
q-G-2.$

We proved that $\{G+1,G+2,\ldots ,q-2\}\subseteq
\mathcal{S}\cup \mathcal{R} . $ Also we showed that the points
$\overline{A}_{t}\ $of $\mathcal{C} \smallsetminus
\mathcal{D}_{G}$ with $t\in \mathcal{S}\cup \mathcal{R}$ are
covered by\ bisecants of $\mathcal{W}_{G}$ either through $Z$
(if $t\in \mathcal{R}$) or\ through $B_{G}$ (if $t\in
\mathcal{S}$). In the other side,\ $\mathcal{C}\smallsetminus
\mathcal{D}_{G}=\{\overline{A} _{t}:t=G+1,G+2,\ldots ,q-2\}\cup
\{A_{\infty },A_{0}\}$ where $\{A_{\infty },A_{0}\}\subset \ell
_{1}$. So, all points of $\mathcal{C}\smallsetminus
\mathcal{D}_{G}$ are covered$.$ $\hfill \qed$
\end{pf}

\begin{definition}
\label{def5.2_Gq}Let $q$ be a \emph{prime power. }Let
$q\not\equiv 3$ $( \bmod~4).$\emph{\ }For integer
$\overline{G},$ let
\begin{equation*}
\mathcal{Z}_{\overline{G}}=\{Z\}\cup \{\overline{A}_{d}:d=0,1,2,\ldots ,
\overline{G}\}.
\end{equation*}
We call \emph{critical value} \emph{of }$\overline{G}$ and
denote by $ \overline{G}_{q}$ the \emph{smallest} value of
$\overline{G}$ such that all points $(1,a,b)$ with $a\in
F_{q}^{\ast },b\in F_{q},b\neq a^{2},$ and all points $(0,1,b)$
with $b\in F_{q}^{\ast },$ lie on bisecants of $\mathcal{Z}
_{\overline{G}}.$
\end{definition}

\begin{thm}
Let $q\geq 32$ be a \emph{prime power. }Let $q\not\equiv 3$
$(\bmod~4).$ If $\overline{G}_{q}\leq \left\lceil
\frac{1}{2}(q-3)\right\rceil $ and
\begin{equation*}
\max \{\overline{G}_{q},\left\lfloor \frac{q-1}{3}\right\rfloor \}\leq G\leq
\left\lceil \frac{q-3}{2}\right\rceil ,
\end{equation*}
then the arc $\mathcal{W}_{G}$ of Construction B is complete.
\end{thm}

\begin{pf}
We use Theorem \ref{th5.2_l0||C-VM_onBisecKM} and Definition
\ref{def5.2_Gq}. $\hfill \qed$
\end{pf}

\noindent In this subsection, we put $q\geq 32$ as we checked
by computer that $\left\lceil \frac{1}{2}(q-3)\right\rceil
<\overline{G}_{q}$ if $q\leq 29.$

\begin{corollary}
\label{cor5.2_family}Let $q\not\equiv 3$ $(\bmod~4)$ be a
\emph{prime power. }If $\overline{G}_{q}\leq \left\lceil
\frac{1}{2}(q-3)\right\rceil $ then Construction B forms a
family of complete $k$-arcs in $PG(2,q)$ containing arcs of all
sizes $k$ in the region
\begin{equation*}
\max \{\overline{G}_{q},\left\lfloor \frac{q-1}{3}\right\rfloor \}+3\leq
k\leq \left\lceil \frac{q-3}{2}\right\rceil +3=\left\{
\begin{array}{cl}
\frac{1}{2}(q+3) & \text{if }q\text{ odd} \\
\frac{1}{2}(q+4) & \text{if }q\text{ even}
\end{array}
\right. .
\end{equation*}
If $\overline{G}_{q}\leq \left\lfloor
\frac{1}{3}(q-1)\right\rfloor $ then size of the smallest
complete arc of the family is $\left\lfloor \frac{1}{3}
(q+8)\right\rfloor .$
\end{corollary}

By computer search using Lemmas \ref{lem5.1_01b_i+j},
\ref{lem5.2_(0,1,b)}, \ref{lem5.2_01b_i+j} we obtained the
following theorem.

\begin{thm}
\label{th5.2_byTable6}Let $q\not\equiv 3$ $(\bmod~4)$ be a
\emph{ prime power. \smallskip }Let $\emph{\ }\overline{G}_{q}$
be given by Definition \emph{\ref{def5.2_Gq}}.\emph{\ }We
introduce $d_{q}$ and $\delta _{q}$ as follows:
$\overline{G}_{q}=d_{q}\sqrt{q}\ln ^{0.95}q,$ $\delta
_{q}=\left\lfloor \frac{1}{3}(q-1)\right\rfloor
-\overline{G}_{q}.$ Then the following holds.
\begin{eqnarray}
\emph{(i)\ }&&\left\lfloor \frac{q-1}{3}\right\rfloor <\overline{G}_{q}\leq
\left\lceil \frac{q-3}{2}\right\rceil \text{ if }32\leq q\leq 81\text{ and }
q=97,101,125. \notag\\
\text{\emph{(ii) }}&&\overline{G}_{q} \leq \left\lfloor \frac{q-1}{3}
\right\rfloor ,~\text{ if }128\leq
q\leq 1367,~
2003\leq q\leq2063,\,\,q =89,109,113,121.\notag\\
\text{\emph{(iii) }}&&\overline{G}_{q}<0.92\sqrt{q}\ln ^{0.95}q~\text{ if }
32\leq q\leq 1367,\,2003\leq q\leq2063.  \label{eq5.2_ln^0.9}\\
\text{\emph{(iv) }}&&
\begin{array}{cll}
1\leq \delta _{q}\leq 93, & 0.66<d_{q}<0.89, & \text{if }128\leq q\leq 593;
\\
89\leq \delta _{q}\leq 196, & 0.71<d_{q}<0.92, & \text{if }601\leq q\leq
1049; \\
183\leq \delta _{q}\leq 267, & 0.68<d_{q}<0.88, & \text{if }1061\leq q\leq
1361;\\
401\leq \delta _{q}\leq 464, & 0.70<d_{q}<0.88, & \text{if }2003\leq q\leq2063.
\end{array} \notag
\end{eqnarray}
\end{thm}

For situations $\left\lfloor \frac{1}{3}(q-1)\right\rfloor
<\overline{G} _{q}, $ the values of $\overline{G}_{q}$ are
given in Table 5.

\begin{thm}\label{th5.2_result}
Let $q\not\equiv 3$ $(\bmod~4)$ be a \emph{prime power. }Let
$128\leq q\leq 1367$, $2003\le q \le2063$, or
$q=89,109,113,121.$ Then Construction B forms a family of
complete $k$-arcs in $PG(2,q)$ containing arcs of all sizes $k$
in the region
\begin{equation*}
\left\lfloor \frac{q+8}{3}\right\rfloor \leq k\leq \left\{
\begin{array}{cl}
\frac{1}{2}(q+3) & \text{if }q\text{ odd} \\
\frac{1}{2}(q+4) & \text{if }q\text{ even}
\end{array}
\right. .
\end{equation*}
\end{thm}

\begin{pf}
We use Corollary \ref{cor5.2_family} and Theorem
\ref{th5.2_byTable6}(ii).$ \hfill \qed$
\end{pf}

\subsection{Arcs with three points outside a conic}

Throughout this subsection, $q\geq 27$ is a \emph{prime power
}and also $ q\equiv 3$ $(\bmod~4).$

Let $J$ be an integer in the region
\begin{equation}
\frac{q-3}{4}\leq J\leq \frac{q-3}{2}.  \label{eq5.3_J}
\end{equation}
Notations $\mathcal{D}_{J},$ $B_{J},$ and $\beta _{J}$ are
taken from (\ref {eq5.2_DG}) and (\ref{eq5.2_Z_BG}) with
substitution $G$ by $J.$ Using (\ref {eq5.3_J}), it is easy to
see that Corollary \ref{cor5.2_(0,1,0)}(i),(iv), Theorem
\ref{th5.2_is-arc} and their proofs hold for $\mathcal{D}_{J},$
$ B_{J},$ and $\beta _{J}$ as well as for
$\mathcal{D}_{G},B_{G},$ and $\beta _{G}.$

\noindent \textbf{Construction C. }Let $q\equiv 3$ $(\bmod~4)$
be a \emph{prime power.}$.$\emph{\ }Let $P,J,$
$\mathcal{D}_{J},$ $Z,$ and $B_{J}$ be given by
(\ref{eq5.1_PT_M}),(\ref{eq5.3_J}),(\ref{eq5.2_DG}), and (\ref
{eq5.2_Z_BG}). We construct a point $(J+4)$-set
$\mathcal{E}_{J}$ in $ PG(2,q) $ as follows:
\begin{equation*}
\mathcal{E}_{J}=\mathcal{D}_{J}\cup \{P,Z,B_{J}\}.
\end{equation*}

\begin{thm}
The $(J+4)$-set $\mathcal{E}_{J}$ of Construction C is an arc.
\end{thm}

\begin{pf}
The set $\mathcal{D}_{J}\cup \{Z,B_{J}\}$ is an arc due to
Theorem \ref {th5.2_is-arc}. By Corollary
\ref{cor5.2_(0,1,0)}(i),(iv), the point $P$ does not lie on
bisecants of $\mathcal{D}_{J}$ and $\mathcal{D}_{J}\cup
\{Z,B_{J}\}$. Finally, $P,Z,B_{J}$ are not collinear.$\hfill
\qed$
\end{pf}

\begin{thm}
\label{th5.3_l0||C-DJ_onBisecEJ}Let $q\equiv 3$ $(\bmod~4)$ be
a \emph{prime power. }Let $J$ be given
by~\emph{(\ref{eq5.3_J})}. Then all points of $\ell _{1}\cup
\mathcal{C}\smallsetminus \mathcal{D}_{J}$ lie on bisecants of
the arc $\mathcal{E}_{J}$ of Construction~C.
\end{thm}

\begin{pf}
All points of $\ell _{1}$ are covered as two points $Z$ and
$B_{J}$ of this line belong to $\mathcal{E}_{J}.$

Throughout this proof, $\mathcal{R},\mathcal{S},$ and
$\mathcal{T}$ are sets of integers modulo $q-1.$ It can be said
that $\mathcal{R},\mathcal{S},$ and $\mathcal{T}$ are sets of
indexes of powers of $\xi .$ We act similarly to the proof of
Theorem \ref{th5.2_l0||C-VM_onBisecKM}.

Let $\mathcal{R}=\{-J,-(J-1),\ldots
,-1\}=\{q-1-J,q-1-(J-1),\ldots ,q-2\}.$ By\
Lemma~\ref{lem5.2_01b_i+j}, a point $\overline{A}_{t}$ of
$\mathcal{C} \smallsetminus \mathcal{D}_{J}$ with $t\in
\mathcal{R}$ lies on the
\smallskip bisecant of $\mathcal{E}_{J}$ through $\overline{A}_{-t}$ and $Z$
where $-t\in \{J,J-1,\ldots ,1\},$ $\overline{A}_{-t}\in
\mathcal{D}_{J}.$

Let $\mathcal{S}=\{\beta _{J}-J,\beta _{J}-(J-1),\ldots ,\beta
_{J}-1,\beta _{J}\}.$ By Lemma \ref{lem5.2_01b_i+j},\ a point
$\overline{A}_{t}\ $of $ \mathcal{C}\smallsetminus
\mathcal{D}_{J}$ with $t\in \mathcal{S}$ lies on the bisecant
$B_{J}\overline{A}_{\beta _{J}-t}$ where $\beta _{J}-t\in
\{J,J-1,J-2,\ldots ,1,0\}$, $\overline{A}_{\beta _{J}-t}\in
\mathcal{D}_{J}.$

Let $\mathcal{T}=\{\frac{1}{2}(q-1),\frac{1}{2}(q-1)+1,\ldots
,\frac{1}{2} (q-1)+J\}.$ By\
(\ref{eq5.2_gamma_-1_1}),(\ref{eq5.2_b=xi^i+xi^j}), points $
P,\overline{A}_{t},$ and $\overline{A}_{t+(q-1)/2}$ are
collinear. Therefore, a point $\overline{A}_{t}$ of
$\mathcal{C}\smallsetminus \mathcal{ D}_{J}$ with $t\in
\mathcal{T}$ lies on the \smallskip bisecant $P\overline{A
}_{t+(q-1)/2}$ where $t+\frac{1}{2}(q-1)\in \{q-1,q,\ldots
,q-1+J\}=\{0,1,\ldots ,J\},$ $\overline{A}_{t+(q-1)/2}\in
\mathcal{D}_{J}.$

As $\beta _{J}=2J,$ we have $\mathcal{S}=\{J,J+1,\ldots ,2J\}$
where $2J\geq \frac{1}{2}(q-1)-1,$ see (\ref{eq5.3_J}). Also,
by (\ref{eq5.3_J}), $\frac{1 }{2}(q-1)+J\geq \frac{1}{4}(3q-5)$
while $q-1-J\leq \frac{1}{4}(3q-5)+1.$

We proved that $\{J+1,J+2,\ldots ,q-2\}\subseteq
\mathcal{S}\cup \mathcal{R} \cup \mathcal{T}.$ Also we showed
that the points $\overline{A}_{t}\ $of $
\mathcal{C}\smallsetminus \mathcal{D}_{J}$ with $t\in
\mathcal{S}\cup \mathcal{R}\cup \mathcal{T}$ are covered
by\smallskip\ bisecants of $ \mathcal{E}_{J}$ either through
$Z$ (if $t\in \mathcal{R}$) or\smallskip\ through $B_{J}$ (if
$t\in \mathcal{S}$) or, finally, through $P$ (if $t\in
\mathcal{T}$). In the other side,\ $\mathcal{C}\smallsetminus
\mathcal{D} _{J}=\{\overline{A}_{t}:t=J+1,J+2,\ldots ,q-2\}\cup
\{A_{\infty },A_{0}\}$ where $\{A_{\infty },A_{0}\}\subset \ell
_{1}$. So, all points of $\mathcal{C }\smallsetminus
\mathcal{D}_{J}$ are covered$.$ $\hfill \qed$
\end{pf}

\begin{definition}
\label{def5.3_Jq}Let $q\equiv 3$ $(\bmod~4)$ be a \emph{prime
power. } For integer $\overline{J},$ let
\begin{equation*}
\mathcal{Q}_{\overline{J}}=\{P,Z\}\cup \{\overline{A}_{d}:d=0,1,2,\ldots ,
\overline{J}\}.
\end{equation*}
We call \emph{critical value} \emph{of }$\overline{J}$ and
denote by $ \overline{J}_{q}$ the \emph{smallest} value of
$\overline{J}$ such that all points $(1,a,b)$ with $a\in
F_{q}^{\ast },b\in F_{q},b\neq a^{2},$ and all points $(0,1,b)$
with $b\in F_{q}^{\ast },$ lie on bisecants of $\mathcal{Q}
_{\overline{J}}.$
\end{definition}

\begin{thm}
Let $q\geq 27$ be a \emph{prime power. }Let $q\equiv 3$
$(\bmod~4).$ If $\overline{J}_{q}\leq \frac{1}{2}(q-3)$ and
\begin{equation*}
\max \{\overline{J}_{q},\frac{q-3}{4}\}\leq J\leq \frac{q-3}{2},
\end{equation*}
then the arc $\mathcal{E}_{J}$ of Construction C is complete.
\end{thm}

\begin{pf}
We use Theorem \ref{th5.3_l0||C-DJ_onBisecEJ} and Definition
\ref{def5.3_Jq}. $\hfill \qed$
\end{pf}

In this subsection, we put $q\geq 27$ as we checked by computer
that $\frac{1 }{2}(q-3)<\overline{J}_{q}$ if $q\leq 23.$

\begin{corollary}
\label{cor5.3_family}Let $q\equiv 3$ $(\bmod~4)$ be a
\emph{prime power. }If $\overline{J}_{q}\leq \frac{1}{2}(q-3)$
then Construction C forms a family of complete $k$-arcs in
$PG(2,q)$ containing arcs of all sizes $k$ in the region
\begin{equation*}
\max \{\overline{J}_{q},\frac{q-3}{4}\}+4\leq k\leq \frac{q+5}{2}.
\end{equation*}
If $\overline{J}_{q}\leq \frac{1}{4}(q-3)$ then cardinality of
this family is equal to $\frac{1}{4}(q-3)$ and size of the
smallest complete arc of the family is $\frac{1}{4}(q+13).$
\end{corollary}

By computer search using Lemmas \ref{lem5.1_01b_i+j},
\ref{lem5.2_(0,1,b)}, \ref{lem5.2_01b_i+j} we obtained the
following theorem.

\begin{thm}
\label{th5.3_byComput}Let $q\equiv 3$ $(\bmod~4)$ be a
\emph{prime power. \smallskip }Let $\emph{\ }\overline{J}_{q}$
be given by Definition~\emph{\ref{def5.3_Jq}}.\emph{\ }We
introduce $r_{q}$ and $\theta _{q}$ as follows:
$\overline{J}_{q}=r_{q}\sqrt{q}\ln ^{0.95}q,$ $\theta
_{q}=\frac{1}{ 4}(q-3)-\overline{J}_{q}.$ Then it holds that
\begin{eqnarray}
\emph{(i)\ }&&\frac{q-3}{4} <\overline{J}_{q}\leq \frac{q-3}{2}\text{ if }
27\leq q\leq 191\text{ and }q=211,223,343;\notag \\
\text{\emph{(ii) }}&&\overline{J}_{q} \leq \frac{q-3}{4},\text{ if }347\leq q\leq 1367,\,
2003\le q\le2063,\notag\\
&&q =199,227,239,243,251,263,271,283,307,311,331;\notag\\
\text{\emph{(iii) }}&&\overline{J}_{q}<0.98\sqrt{q}\ln ^{0.95}q\text{ if }
27\leq q\leq 1367,\,2003\le q\le2063;  \label{eq5.3_098ln^095}\\
\text{\emph{(iv) }}&&
\begin{array}{cll}
3\leq \theta _{q}\leq 44, & 0.67<r_{q}<0.86, & \text{if }347\leq q\leq 599;
\\
18\leq \theta _{q}\leq 111, & 0.67<r_{q}<0.94, & \text{if }607\leq q\leq 991;
\\
89\leq \theta _{q}\leq 167, & 0.67<r_{q}<0.84, & \text{if }1019\leq q\leq
1367;\\
261\leq \theta _{q}\leq 294, & 0.67<r_{q}<0.80, & \text{if }2003\le q\le 2063.
\end{array} \notag
\end{eqnarray}
\end{thm}

For situations $\frac{1}{4}(q-3)<\overline{J}_{q},$ the values
of $\overline{ J}_{q}$ are given in Table 5.

\begin{thm}\label{th5.3_result}
Let $q\equiv 3$ $(\bmod~4)$ be a \emph{prime power. }Let
$347\leq q\leq 1367$ or
$q=199,227,239,243,251,263,$ $271,283,307,311,331.$ Then
Construction C forms a family of complete $k$-arcs in $PG(2,q)$
containing arcs of all sizes $k$ in the region
\begin{equation*}
\frac{q+13}{4}\leq k\leq \frac{q+5}{2}.
\end{equation*}
\end{thm}

\begin{pf}
We use Corollary \ref{cor5.3_family} and Theorem
\ref{th5.3_byComput}(ii).$ \hfill \qed$
\end{pf}

Basing on Theorems \ref{th5.1_byTable1}, \ref{th5.2_byTable6},
\ref{th5.3_byComput} and taking into account that $\sqrt{q}\ln
^{0.9}q<\sqrt{q}\ln q$, $\sqrt{q}\ln ^{0.95}q<\sqrt{q}\ln q$,
we conjecture the following, cf. Conjecture
\ref{conj1_construc}.

\begin{conjecture}\label{conj5_constr}
Let $\overline{H}_{q},$ $\overline{G}_{q}$, $\overline{J}_{q}$
be given by Definitions \emph{\ref{def5.1_overlineH}},
\emph{\ref {def5.2_Gq}}, \emph{\ref{def5.3_Jq}}. Let for
$\overline{H}_{q},$ $q$ be prime while for $\overline{G}_{q}$
and $\overline{J}_{q}$ it holds that $q$ is a prime power.
Finally, let $q\not\equiv 3$ $(\bmod~4)$ for $\overline{G}_{q}$
and $q\equiv 3$ $(\bmod~4)$ for $\overline{J}_{q}$. Then the
following holds.
\begin{equation*}
\overline{H}_{q}\leq \left\lfloor \frac{q-1}{3}
\right\rfloor\text{ if }
q\geq 109;~\overline{G}_{q}\leq \left\lfloor \frac{q-1}{3}
\right\rfloor\text{ if }q\geq 128;~
\overline{J}_{q}\leq \frac{q-3}{4}\text{ if }q\geq 347.
\end{equation*}
\begin{equation}
\overline{H}_{q}<\sqrt{q}\ln q\text{ if }q\geq 19;~
\overline{G}_{q}<\sqrt{q}\ln q\text{ if }q\geq 32;~
\overline{J}_{q}<\sqrt{q}\ln q\text{ if }q\geq 27.
\label{eq5.3_conject}
\end{equation}
\end{conjecture}

\begin{rmk}
It is interesting to compare the relations
(\ref{eq5.1_ln^0.9}),(\ref{eq5.2_ln^0.9}),(\ref
{eq5.3_098ln^095}),(\ref{eq5.3_conject}) with Theorems
\ref{th1_ln0.75}, \ref{th4_ln0.75}, \ref{th4_ln0.741} and
Conjecture \ref{conj1_ln0.75} and to compare also computer
results providing Theorems \ref {th5.1_byTable1},
\ref{th5.2_byTable6}, \ref{th5.3_byComput} with Tables 1 and 2.
One can see that the upper estimates of $t_{2}(2,q),$
$\overline{H} _{q},$ $\overline{G}_{q},$ and $\overline{J}_{q}$
have the same structure and the values of
$\overline{t}_{2}(2,q),$ $\overline{H}_{q},$
$\overline{G}_{q},$ and $ \overline{J}_{q}$ have a close order.
This seems to be natural as \emph{ almost all }points of
$PG(2,q)$ lie on bisecants of $\mathcal{P}_{\overline{H
}_{q}},$ $\mathcal{Z}_{\overline{G}_{q}},$ and
$\mathcal{Q}_{\overline{J} _{q}},$ see Definitions
\ref{def5.1_overlineH}, \ref{def5.2_Gq} and \ref {def5.3_Jq}.
\end{rmk}

\begin{rmk}\label{rem5_starting-objects}
The complete arcs of Constructions A, B, C can be used
as starting objects in inductive constructions.  For example, for even
$q$, arcs of Construction~B can be used in constructions of
\cite[Ths 1.1,3.14-3.17,4.6-4.8]{DGMP-JCD}. In that
way, one can generate infinite sets of families of complete caps in
projective spaces $PG(v,2^{n})$ of growing dimensions $v$. For every $v$,
constructions of \cite{DGMP-JCD} can obtain a complete cap from every complete arc
of Construction B. Also, it can be shown that in Constructions A,
B, C all points not on conic are \emph{external}. So, the arcs of
Constructions A and C for $q\equiv 3~(\bmod~4)$ can be used as
starting objects in constructions of \cite{KorchPace}, see
\cite[Th.\thinspace 23]{KorchPace}. Thereby, infinite families of large
complete arcs in $PG(2,q^{n})$ with growing $n$ can be obtained.
\end{rmk}

\section{On the spectrum of possible sizes of complete arcs in $PG(2,q)$
\label{sec6_spectrum}}

The main known results on the spectrum of possible sizes of complete arcs in
$PG(2,q)$ are given in Introduction with the corresponding references.
Taking into account the results cited in Introduction, we denote
\begin{equation*}
M_{q}=\left\{
\begin{array}{ll}
\frac{1}{2}(q+4) & \text{for even }q\smallskip \\
\frac{1}{2}(q+7) & \text{for odd }q\text{ included to (\ref{eq1+(q+7)/2})}
\smallskip \\
\frac{1}{2}(q+5) & \text{for odd }q\text{ not included to (\ref{eq1+(q+7)/2}
) }
\end{array}
\right. .
\end{equation*}
We suppose that the smallest known sizes
$\overline{t}_{2}(2,q)$ are given in Tables 1-4 of this paper.

\begin{thm}
\label{th5_spectrum} In $PG(2,q)$ with $25\leq q\leq 349 $, and
$q=1013,$ $2003$, there are complete $k$-arcs of \textbf{all}
the sizes in the region $$ \overline{t}_{2}(2,q)\leq k\leq
M_{q}.$$
\end{thm}

\begin{pf}
For $25\le q\le 167$ the assertion of the theorem follows from
\cite[Tab.\thinspace2]{DFMP-JG2005} and
\cite[Tab.\thinspace2]{DFMP-JG2009}. For $169\le q\le 349$ and
$q=1013,2003$, we used Constructions A, B, C of Section
\ref{sec5_constr}. The sizes not following from the
constructions are obtained in this work by the randomized
greedy algorithms.
\end{pf}

An experience obtained in computer search for the proof of
Theorem \ref {th5_spectrum} allows us to do Conjecture
\ref{conj5}. Here we took into account sizes that can be got by
Constructions A, B, C and the remark on sizes close to
$\overline{t}_{2}(2,q)$ in the end of Section
\ref{sec2_computSearch}. Note also that the rest of sizes  for
all $169\le q\le 349$ and $q=1013,2003$ was relatively easy
obtained by the greedy algorithms with point subset of a conic
taken as the starting set $S_{0}$. For this we used
consequently subsets of cardinality approximately
$15\%,20\%,25\%,30\%$ of the conic cardinality $q+1$.
\begin{conjecture}
\label{conj5} Let $353\le q\le 4561$, $ q\in T_{2}\cup T_{3}$,
be  a prime power. Then in $PG(2,q)$ there are complete
$k$-arcs of \textbf{all} the sizes in the region $
\overline{t}_{2}(2,q)\le k\le M_{q}$. Moreover, complete
$k$-arcs with $ \overline{t}_{2}(2,q)\le k\le \frac{1}{2}(q+5)$
can be obtained either by Constructions A, B, C or by the
randomized greedy algorithms.
\end{conjecture}

\end{document}